%% file: main.tex
\documentclass[onefignum,onetabnum]{siamonline220329}

\input{main_shared}

\ifpdf
\hypersetup{
  pdftitle={Decomposition of Boolean networks: An approach to modularity of biological systems},
  pdfauthor={Claus Kadelka, Reinhard Laubenbacher, David Murrugarra, Alan Veliz-Cuba, Matthew Wheeler}
}
\fi

\begin{document}

\maketitle

% REQUIRED
\begin{abstract}
This paper presents the foundation for a decomposition theory for Boolean networks, a type of discrete dynamical system that has found a wide range of applications in the life sciences, engineering, and physics. Given a Boolean network satisfying certain conditions, there is a unique collection of subnetworks so that the network can be reconstructed from these subnetworks by an extension operation. The main result of the paper is that this structural decomposition induces a corresponding decomposition of the network dynamics. The theory is motivated by the search for a mathematical framework to formalize the hypothesis that biological systems are modular, widely accepted in the life sciences, but not well-defined and well-characterized. As an example of how dynamic modularity could be used for the efficient identification of phenotype control, the control strategies for the network can be found by identifying controls in its modules, one at a time. 
\end{abstract}

% REQUIRED

% REQUIRED

\section{Introduction}

Building complicated structures from simpler building blocks is a widely observed principle in both natural and engineered systems. In molecular systems biology, it is also widely accepted, even though there has not emerged a clear definition of what constitutes a simple building block, or module. Consequently, it is not clear how the modular structure of a system can be identified, why it is advantageous to an organism to be composed of modular systems, and how we could take advantage of modularity to advance our understanding of molecular systems; see, e.g., \cite{Hartwell1999, hernandez2022effects, Lorenz2011}. In the (graph-theoretic) network representation of molecular systems, such as gene regulatory networks or protein-protein interaction networks, a module refers to a ``highly" connected region of the graph that is ``sparsely" connected to the rest of the graph, otherwise known as communities in graphs. Graph search algorithms that depend on the choice of parameters are used to define modules. Similar approaches are used for identifying modules in co-expression networks based on clustering of transcriptomics data. 

A further limitation of this approach to modularity is that it focuses entirely on static representations of gene regulatory networks and other systems. However, living organisms are dynamic, and need to be modeled and understood as dynamical systems. Modularity should have an instantiation as a dynamic feature. The most common types of models employed for this purpose are systems of ordinary differential equations and discrete models such as Boolean networks and their generalizations, providing the basis for a study of dynamic modularity. In recent years, there have been an increasing number of papers that take this point of view. The authors of \cite{jimenez2017modularity} argue that dynamic modularity may be independent of structural modularity and identify examples of multifunctional circuits in gene regulatory networks that they consider dynamically modular but without any structural modularity underpinning this capability. A similar argument is made in \cite{jaeger2019modularity} by analyzing a small gene regulatory network example. For another example of a similar approach see \cite{Deritei2016}. 

The literature on how modularity might have evolved and why it might be useful as an organizational principle cites as the most common reasons robustness, the ability to rapidly respond to changing environmental conditions, and efficiency in the control of response to perturbations. An interesting hypothesis has been put forward in \cite{Lorenz2011}, namely that a modular organization of biological structure can be viewed as a symmetry-breaking phase transition, with modularity as the order parameter. 

This literature makes clear that research on the topic of modularity in molecular systems, both structural and dynamic, would be greatly advanced by clear definitions of the concept of module, both structural and dynamic. This would in particular help to decide whether and how structural and dynamic modularity are related, and it would provide a basis on which to distinguish between modularity and multistationarity of a dynamic regulatory network. To be of practical use, such a theory should include algorithms to decompose a dynamic network into structural and/or dynamic modules. At the same time, it would be of great practical value, for instance for synthetic biology, to understand how systems can be composed from modules that have specific dynamic properties.

The desire for algorithms has led us to look for guidance to mathematics, as a complement to biology. Many of the important mathematical theories in a wide range of fields are decomposition theories (e.g., algebra, topology, or geometry). Choosing the widely-used modeling framework of Boolean networks, we asked whether it is possible to identify meaningful concepts of modularity that, ideally, link both the structural and dynamic aspects. The concept of module we propose is structural, in terms of strongly connected components of the wiring diagram. This leads to an unambiguous decomposition into non-overlapping modules (Section~\ref{sec:dec_mod}). The surprising consequence is that this structural property also induces a decomposition of the network dynamics (Section~\ref{sec:dec_dyn}). It also implies that the modules are organized in a hierarchical fashion, which has interesting implications for the control of such networks (Section~\ref{sec:control}). This provides evidence for the hypothesized advantages of a modular organization for an organism. 

The main question to be answered at this point, unfortunately beyond the scope of this paper, is whether there is biological evidence that our concept of module and our dynamic decomposition theory does in fact reflect biological reality. As a preliminary study, we have analyzed all 130 models in the database \cite{kadelka2020meta}, with up to 320 nodes. The vast majority of them satisfy the criteria that define a module, qualifying them as ``modular'' in our theory. This raises the question why not more of these models have more than one module. One possible interpretation of this observation is that modules should be subnetworks that carry out key control functions in a cell. It would therefore not be surprising if there was a selection bias among systems biologists to focus their attention on such modules, given that larger networks are still challenging to build and analyze. A good starting point for a search for modules might be genome-scale co-expression networks that tend to show a larger number of modules (defined using graph connectivity measures). 

For the purpose of this article, we will focus on the class of Boolean networks as a modeling paradigm. 
Recall that a Boolean network $F$ on variables $x_1, \ldots , x_n$ can be viewed as a function on binary strings of length $n$, which is described coordinate-wise by $n$ Boolean rules $f_i$ in the $x_i$ (or a subset thereof). Two directed graphs can be associated to $F$: one is the \textit{wiring diagram} (also known as dependency graph) with nodes corresponding to the $x_i$. There is a directed edge from $x_i$ to $x_j$ if $f_j$ depends on $x_i$. The \textit{state space} of $F$ has as nodes the $2^n$ binary strings. There is a directed edge from $\mathbf{u}$ to $\mathbf{v}$ if $F({\mathbf u}) = \mathbf v$. Each connected component of the state space gives an attractor basin of $F$, a directed loop as a limit cycle, with trees feeding into the nodes in the loop. 

Here we develop a concept of module for Boolean networks and a decomposition theory for Boolean networks that have special functional properties. We initially focused Boolean networks governed by so-called nested canalizing functions, a concept first introduced in \cite{KauffmanNCF2003}, due to its implications for network dynamics. We and others have studied this class of Boolean functions and their networks extensively (see, e.g., \cite{jarrah2007nested, he2016stratification, kadelka2017influence, li2013boolean}). This class is restrictive enough and admits sufficiently rich mathematical structure to enable a decomposition theory, yet it is sufficiently general to capture almost all regulatory relationships that appear in published Boolean network models of GRNs \cite{kadelka2020meta}. Nested canalizing functions, in turn, are a special case of fanout-free Boolean functions. These have the property that they can be written in logical form so that each variable only occurs once. They form a subclass of the class of monotone functions, for which each variable is either an inhibitor or an activator, a biologically meaningful property. The condition of being fanout-free is sufficient for the decompositions described here.

%We and others have studied so-called nested canalizing functions, demonstrated to be common in gene regulatory network (GRN) models~\cite{kauffman2004genetic,li2013boolean,kadelka2017influence,jarrah2007nested}. They, in turn, are a special case of fanout-free Boolean functions. These have the property that they can be written in logical form so that each variable only occurs once. They form a subclass of the class of monotone functions, for which each variable is either an inhibitor or an activator, a biologically meaningful property.  

%We initially focused on so-called nested canalizing Boolean networks, a concept first introduced in \cite{KauffmanNCF2003}, due to its implications for network dynamics. We and others have studied this class of Boolean functions and their networks extensively (see, e.g., \cite{he2016stratification, kadelka2017influence, li2013boolean} ). This class is restrictive enough and admits sufficiently rich mathematical structure to enable a decomposition theory, yet it is sufficiently general to capture almost all regulatory relationships that appear in published Boolean network models of GRNs \cite{kadelka2020meta}. However, the weaker condition of being fanout-free is sufficient. 

Our concept of modularity is structural. A module of a Boolean network is a subnetwork that has a strongly connected wiring diagram. The strongly connected components are disjoint. Two components can be connected by edges going in one direction only. However, what might be considered our main result, this structural decomposition induces a dynamic decomposition that allows the reconstruction of the dynamics of the entire network from the dynamics of its modules in an algorithmic way.  Conversely, we show how to construct networks by combining Boolean networks into a larger network that has them as modules. These two constructions, decomposition and combination, are inverse to each other for fanout-free and linear networks.

One rationale proposed for the modularity principle is that modules are key control units that are largely evolutionarily preserved and exist across species, whereas the assembly of modules provides evolutionary freedom. In order to provide an estimate for this evolutionary freedom we parametrized the different ways in which to assemble a given set of modules into networks, in Section~\ref{sec:enumeration}. Finally, in Section~\ref{sec:control}, we show how modularity has advantages for the purpose of control. We show that control targets can be identified one module at a time, taking into account the downstream effect of controls in upstream components. This is demonstrated with an example of a regulatory network important in cancer.

\section{Background}
In this section, we review some standard definitions and introduce the concepts of \emph{canalization} and \emph{fanout-freeness}. Throughout the paper, let $\mathbb F_2$ be the binary field with elements 0 and 1. %Further, let $\oplus$ denote addition modulo $2$ when used in a polynomial, and the ``exclusive or" (XOR) function when used in a Boolean logical expression. \claus{should we even bother using $\oplus$? thus far, we are inconsistent in its use}

\subsection{Boolean functions}

%\subsection{Fanout-Free Functions}
The main class of Boolean functions that we consider in this paper will be fanout-free functions. 
\begin{definition}\label{fanout-free} A Boolean function is {\emph fanout-free} if it can be written in logical format (using $\wedge, \vee, \neg$ only) so that each variable appears only once.
\end{definition}

\begin{example}\label{ex:first}
The Boolean functions 
$f(x_1,x_2,x_3,x_4)=x_1\wedge (\neg x_2\vee (x_3 \wedge x_4))$ as well as $g(x_1,x_2,x_3,x_4)=(x_1\wedge \neg x_2)\vee (x_3 \wedge x_4)$ are fanout-free. The threshold function 
\[h(x_1,x_2,x_3)= (x_1\wedge x_2)\vee(x_1\wedge x_3)\vee(x_2\wedge x_3) = \begin{cases} 1 & \text{if }\sum_i x_i > 1,\\0 & \text{otherwise,}\end{cases}\] is not fanout-free.
%Note that $g$ is not nested canalizing. From the representation in Theorem~\ref{thm:he}, it is straight-forward to see that any nested canalizing function is fanout-free.
\end{example}

The main reason for considering fanout-free functions is that they come with a natural choice of restriction, which we will discuss in more detail in Section~\ref{section:restriction}.  Common examples of fanout-free functions include conjunctive, disjunctive, AND-NOT functions, or more generally, nested canalizing functions.

\begin{definition}\label{def_essential}
A Boolean function $f:\F^n \to \F$ \emph{depends} on the variable $x_i$ if there exists an $\vecx \in \FF_2^n$ such that 
$$f(\vecx)\neq f(\vecx + e_i),$$
where $e_i$ is the $i$th unit vector. %In that case, we also say $f$ \emph{depends} on $x_i$.
\end{definition}

% \begin{example}
% The  is not an example of a fanout-free network as at least one of the variables appears twice
% \end{example}

% \begin{definition}\label{def_increasing}
% A Boolean function $f:\F^n \to \F$ is \emph{increasing} (\emph{decreasing}) in the variable $x_i$ if it depends on $x_i$ and for all $\vecx \in \FF_2^n$ we have
% $$f(x_1,\ldots,x_{i-1},0,x_{i+1},\ldots,x_n) \leq (\geq) f(x_1,\ldots,x_{i-1},1,x_{i+1},\ldots,x_n).$$
% \end{definition}

%\begin{definition}\label{def_canalizing}
%A Boolean function $f:\F^n \to \F$ is \emph{canalizing} if there exists a variable $x_i$, $a,\,b\in\{0,\,1\}$ and a Boolean function $g(x_1,\ldots,x_{i-1},x_{i+1},\ldots,x_n)\not\equiv b$  such that
%$$f(x_1,x_2,...,x_n)= \begin{cases}
%b,& \ \text{if}\ x_i=a\\
%g(x_1,x_2,...,x_{i-1},x_{i+1},...,x_n),& \ \text{if}\ x_i\neq a
%\end{cases}$$
%In that case, we say that $x_i$ \emph{canalizes} $f$ (to $b$) and call $a$ the \emph{canalizing input} (of $x_i$) and $b$ the \emph{canalized output}.
%\end{definition}

\begin{definition}\label{def_nested_canalizing}
A Boolean function $f:\F^n \to \F$ is \emph{nested canalizing} with respect to the permutation $\sigma \in \mathcal{S}_n$, inputs $a_1,\ldots,a_n$ and outputs $b_1,\ldots,b_n$, if
\begin{equation*}f(x_{1},\ldots,x_{n})=
\left\{\begin{array}[c]{ll}
b_{1} & x_{\sigma(1)} = a_1,\\
b_{2} & x_{\sigma(1)} \neq a_1, x_{\sigma(2)} = a_2,\\
b_{3} & x_{\sigma(1)} \neq a_1, x_{\sigma(2)} \neq a_2, x_{\sigma(3)} = a_3,\\
\vdots  & \vdots\\
b_{n} & x_{\sigma(1)} \neq a_1,\ldots,x_{\sigma(n-1)}\neq a_{n-1}, x_{\sigma(n)} = a_n,\\
1 + b_{n} & x_{\sigma(1)} \neq a_1,\ldots,x_{\sigma(n-1)}\neq a_{n-1}, x_{\sigma(n)} \neq a_n.
\end{array}\right.\end{equation*}
The last line ensures that $f$ actually depends on all $n$ variables. From now on, we use the acronym NCF for nested canalizing function.
\end{definition}

We restate the following powerful stratification theorem for reference, which provides a unique polynomial form for any Boolean function.

\begin{theorem}[\cite{he2016stratification}]\label{thm:he}
Every Boolean function $f(x_1,\ldots,x_n)\not\equiv 0$ can be uniquely written as 
\begin{equation}\label{eq:matts_theorem}
%\begin{array}[c]{ll}
    f(x_1,\ldots,x_n) = M_1(M_2(\cdots (M_{r-1}(M_rp_C + 1) + 1)\cdots)+ 1)+ q,
%\end{array}
\end{equation}

where each $M_i = \displaystyle\prod_{j=1}^{k_i} (x_{i_j} + a_{i_j})$ is a nonconstant extended monomial, $p_C$ is the \emph{core polynomial} of $f$, and $k = \displaystyle\sum_{i=1}^r k_i$ is the canalizing depth. Each $x_i$ appears in exactly one of $\{M_1,\ldots,M_r,p_C\}$, and the only restrictions are the following ``exceptional cases'':
\begin{enumerate}
    \item If $p_C\equiv 1$ and $r\neq 1$, then $k_r\geq 2$;
    \item If $p_C\equiv 1$ and $r = 1$ and $k_1=1$, then $q=0$.
\end{enumerate}
When $f$ is not canalizing (\textit{i.e.}, when $k=0$), we simply have $p_C = f$.
\end{theorem}

\begin{example}
From the representation in Theorem~\ref{thm:he}, it is straight-forward to see that the Boolean function $f$ from Example~\ref{ex:first} is nested canalizing, and more generally, that any NCF is fanout-free. Further, $g$ and $h$ from Example~\ref{ex:first} are not nested canalizing.
\end{example}

% \begin{remark}~\label{rem:layer_output} Note the following properties of canalization.

% (a) Theorem~\ref{thm:he} shows that any Boolean function has a unique extended monomial form, in which the variables are partitioned into different layers based on their dominance. Any variable that is canalizing (independent of the values of other variables) is in the first layer. Any variable that ``becomes'' canalizing when excluding all variables from the first layer is in the second layer, etc. Variables in any layer will be referred to as \emph{conditionally canalizing}. All remaining variables that never become canalizing are part of the core polynomial. The number of variables that eventually become canalizing, \textit{i.e.}, the number of conditionally canalizing variables, is the canalizing depth of the function. NCFs are exactly those functions where all variables are conditionally canalizing.

% (b) While variables in the same layer may have different canalizing input values, they all share the same canalized output value, \textit{i.e.}, they all canalize a function to the same output. \new{On the other hand, the outputs of two consecutive layers are distinct. Therefore, the number of layers of a $k$-canalizing function expressed as in Definition~\ref{def:kcanalizing} is simply one plus the number of changes in the vector of canalized outputs, $(b_1,b_2,\ldots,b_r)$.}
% \end{remark}

From Equation~\ref{eq:matts_theorem}, we can directly derive an important summary statistic of NCFs.

\begin{definition}\label{def_layers} \emph{(\cite{kadelka2017influence})}
Given an NCF $f(x_1,\ldots,x_n)$ represented as in Equation~\ref{eq:matts_theorem}, we call the extended monomials $M_i$ the \emph{layers} of f and define the \emph{layer structure} as the vector $(k_1,\ldots,k_r)$, which describes the number of variables in each layer. Note that $k_1+\cdots+k_r = n$ and that by exceptional case 1 in Theorem~\ref{thm:he}, $k_r\geq 2$. 
\end{definition}

\begin{example}
The Boolean functions 
\begin{align*}
    f(x_1,x_2,x_3,x_4)&=x_1\wedge (\neg x_2\vee (x_3 \wedge x_4)) = x_1\left[x_2\left[x_3x_4+1\right]+1\right],\\
    g(x_1,x_2,x_3,x_4)&=x_1\wedge (\neg x_2\vee x_3 \vee x_4) = x_1\left[x_2\left(x_3+1\right)\left(x_4+1\right)+1\right]
\end{align*}
are nested canalizing. $f$ consists of three layers with layer structure $(1,1,2)$, while $g$ possesses only two layers and layer structure $(1,3)$.
\end{example}

%\subsection{Fanout-Free Functions}

%\begin{definition}\label{fanout-free} A Boolean function is {\emph fanout-free} if it can be written in logical format (using $\wedge, \vee, \neg$ only) so that each variable appears only once.
%\end{definition}

%\begin{example}
%The Boolean functions 
%$f(x_1,x_2,x_3,x_4)=x_1\wedge (\neg x_2\vee (x_3 \wedge x_4))$ as well as $h(x_1,x_2,x_3,x_4)=(x_1\wedge \neg x_2)\vee (x_3 \wedge x_4)$ are fanout-free. Note that $h$ is not nested canalizing. From the representation in Theorem~\ref{thm:he}, it is straight-forward to see that any nested canalizing function is fanout-free.
%\end{example}

\subsection{Boolean Networks}

\begin{definition}\label{def_local_model}
  A \emph{Boolean network} is an $n$-tuple of \emph{coordinate functions} $F:=(f_1,\dots,f_n)$, where $f_i\colon \FF_2^n\to \FF_2$. Each function $f_i$ uniquely determines a map
\begin{equation*}\label{eqn:F_i}
  F_i\colon \FF_2^n\longto \FF_2^n,\qquad F_i\colon(x_1,\dots,x_n)\longmapsto (x_1,\dots,{f_i(x)},\dots,x_n),
\end{equation*}
where $x=(x_1,\dots,x_n)$. Every Boolean network defines a canonical map, where the functions are synchronously updated:
\[
F\colon \FF_2^n\longto \FF_2^n,\qquad F\colon(x_1,\dots,x_n)\longmapsto(f_1(x),\dots,f_n(x)).
\]
\end{definition}
In this paper, we only consider this canonical map, i.e., we only consider synchronously updated Boolean network models. Frequently, we use the terms \emph{fanout-free network} or \emph{nested canalizing network}, which simply means Boolean networks where all update rules are from the respective family of functions.

% Every function from a finite set to itself can be represented as a directed graph, where the directed edges are of the form $(x,f(x))$. For Boolean network models, this graph is called the phase space, or state space. It provides a visual representation of the time-evolution of the system.

% \begin{definition}\label{def_phase_space}
%   Let $f=(f_1,\dots,f_n)$ be a Boolean network model. The \emph{phase space} of $f$ is the directed graph with vertex set $\FF_2_2^n$ and edge set $\{(x,f(x))\mid x\in \FF_2_2^n\}$. A node $x$ in the phase space is a:
%   \begin{itemize}
%   \item \emph{transient point} if $f^k(x)\neq x$ for all $k>1$,
%   \item \emph{periodic point} if $f^k(x)=x$ for some $k\geq 1$,
%   \item \emph{fixed point} if $f(x)=x$.
%   \end{itemize}
% \end{definition}

\begin{definition}\label{def_wiring_diagram}
The \emph{wiring diagram} of a Boolean network $F=(f_1,\ldots,f_n): \FF_2^n\rightarrow \FF_2^n$ is the directed graph with vertices $x_1,\ldots,x_n$ and an edge from $x_i$ to $x_j$ if $f_j$ depends on $x_i$. 
\end{definition}

% Since we are interested in networks that represent biological systems, each edge in the wiring diagram will often have a sign associated to it. %\claus{Note that we do see a few ``conditional" edges in my database, i.e. edges that are neither, https://arxiv.org/pdf/2009.01216.pdf Figure 1E.}. 
% In that case, an edge from $i$ to $j$ will be positive if $f_j$ is increasing in $x_i$, and the edge will be negative if $f_j$ is decreasing in $x_i$. %Note: A Boolean function is increasing (decreasing) in a variable if switching this variable's input from 0 to 1 will never decrease (increase) the function output, irrespective of the other inputs. 
% Such networks are called unate networks and include those constructed using fanout-free functions.

\begin{figure}
    \centering
    \includegraphics[scale=1.5]{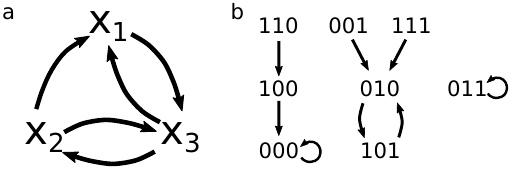}
    \caption{Wiring diagram and state space of the Boolean network in Example~\ref{eg:wd}-\ref{eg:cycle_structure}. 
    (a) The wiring diagram encodes the dependency between variables. %The edges $x_3\rightarrow x_1$ and $x_1\rightarrow x_3$ are negative because $f_1=x_2 \wedge \neg x_3$ is decreasing in $x_3$ and $f_3=\neg x_1\wedge x_2$ is decreasing in $x_1$.
    (b) The state space is a directed graph with edges between all states and their images. This graph therefore encodes all possible trajectories. 
    }
    \label{fig:dynamics_eg}
\end{figure}

\begin{example}\label{eg:wd}
%Consider the Boolean network $F:\FF_2^3\rightarrow \FF_2^3$ given by $F(x)=(\neg x_2, x_3\vee x_2, \neg x_1\wedge x_2)$. Its wiring diagram is shown in Figure~\ref{fig:wd}
Figure~\ref{fig:dynamics_eg}a shows the wiring diagram of the Boolean network $F:\FF_2^3\rightarrow \FF_2^3$ given by $$F(x_1,x_2,x_3)=(x_2 \wedge \neg x_3,x_3,\neg x_1 \wedge x_2).$$
\end{example}

% \begin{figure}
%     \centering
%     \includegraphics[scale=1.1]{wd.pdf}
%     \caption{Wiring diagram of network in Example \ref{eg:wd}. The edges $x_2\rightarrow x_1$ and $x_1\rightarrow x_3$ are negative because $f_1=\neg x_2$ is decreasing with respect to $x_2$ and $f_3=\neg x_1\wedge x_2$ is decreasing with respect to $x_1$. }
%     \label{fig:wd}
% \end{figure}

\subsection{Dynamics of Boolean networks}

There are two ways to describe the dynamics of a Boolean network $F: \FF_2^n\to \FF_2^n$, (i) as trajectories for all $2^n$ possible  initial conditions, or (ii) as a directed graph with nodes in $\FF_2^n$. Although the first description is less compact, it will allow us to formalize the dynamics of coupled networks. 

%\subsubsection{Boolean network dynamics as trajectories}
\begin{definition}
A trajectory of a Boolean network $F:\FF_2^n\rightarrow \FF_2^n$ is a sequence $(x(t))_{t=0}^\infty$ of elements of $\FF_2^n$ such that $x(t+1)=F(x(t))$ for all $t\geq 0$.
\end{definition}

\begin{example}\label{eg:traj}
%Consider the Boolean network $F:\FF_2^n\rightarrow \FF_2^n$ given by 
%$$F(x_1,x_2,x_3)=(x_2 \wedge \neg x_3,x_3,\neg x_1 \wedge x_2).$$
%Its wiring diagram is shown in Figure~\ref{fig:dynamics_eg}a.
For the example network above, $F(x_1,x_2,x_3)=(x_2 \wedge \neg x_3,x_3,\neg x_1 \wedge x_2)$, there are $2^3=8$ possible initial states giving rise to the following trajectories (commas and parenthesis for states are omitted for brevity).
\begin{align*}
     T_1&=(000,000,000,000,\ldots)\\
     T_2&=(001,010,101,010,\ldots)\\
    T_3&=(010,101,010,101,\ldots)\\
    T_4&=(011,011,011,011,\ldots)\\
    T_5&=(100,000,000,000,\ldots)\\
    T_6&=(101,010,101,010,\ldots)\\
    T_7&=(110,100,000,000,\ldots)\\
    T_8&=(111,010,101,010,\ldots)
\end{align*}
% \begin{itemize}
%     \item $T_1=(000,000,000,000,\ldots)$
%     \item $T_2=(001,010,101,010,\ldots)$
%     \item $T_3=(010,101,010,101,\ldots)$ 
%     \item $T_4=(011,011,011,011,\ldots)$
%     \item $T_5=(100,000,000,000,\ldots)$
%     \item $T_6=(101,010,101,010,\ldots)$
%     \item $T_7=(110,100,000,000,\ldots)$ 
%     \item $T_8=(111,010,101,010,\ldots)$
% \end{itemize}
We can see that $T_3$ and $T_6$ are periodic trajectories with period 2. Similarly, $T_1$ and $T_4$ are periodic with period 1. All other trajectories eventually reach one of these 4 states.
\end{example}

When seen as trajectories, $T_3$ and $T_6$ are different, but they can both be encoded by the fact that $F(0,1,0)=(1,0,1)$ and $F(1,0,1)=(0,1,0)$. Similarly, $T_1$ and $T_4$ can be encoded by the equalities $F(0,1,1)=(0,1,1)$, $F(0,0,0)=(0,0,0)$. This alternative, more compact way of encoding the dynamics of a Boolean network is the standard approach, which we formalize next. 

%\subsubsection{Boolean network dynamics as a directed graph with nodes in $\FF_2^n$}

%A compact way of encoding all trajectories of a Boolean network is given in the following definition.

\begin{definition}
The \emph{state space} of a (synchronously updated) Boolean network $F:\FF_2^n\to \FF_2^n$ 
is a directed graph with vertices in $\FF_2^n$ and an edge from $x$ to $y$ if $F(x)=y$.
\end{definition}

\begin{example}\label{eg:state_space}
Figure~\ref{fig:dynamics_eg}b shows the state space of the (synchronously updated) Boolean network from Example~\ref{eg:wd}.
%\[F(x_1,x_2,x_3)=(x_2 \wedge \neg x_3,x_3,\neg x_1 \wedge x_2).\]
\end{example} 

From the state space, one can easily obtain all periodic points, which form the attractors of the network. 

\begin{definition}\label{def_attractors}
The \emph{space of attractors} for a Boolean network is the collection $\mathcal{D}(F)$ of all \emph{minimal} subsets $\mathcal{C}\subseteq \FF_2^n$ satisfying $F(\mathcal{C})=\mathcal{C}$.  
\begin{enumerate}
    \item The subset $\mathcal{D}^1(F)\subset\mathcal{D}(F)$ of sets of exact size 1 consists of all \emph{steady states} (also known as \emph{fixed points}) of $F$.%$(F,X)$.
    \item The subset $\mathcal{D}^r(F)\subset\mathcal{D}(F)$ of sets of exact size $r$ consists of all \emph{cycles} of length $r$ of $F$.%$(F,X)$. 
\end{enumerate}

Equivalently, an \emph{attractor of length $r$} is an ordered set with $r$ elements, $\mathcal{C}=\{c_1,\ldots,c_r\}$, such that $F(c_1)=c_2, F(c_2)=c_3,\ldots, F(c_{r-1})=c_r, F(c_r)=c_1$.
\end{definition}

\begin{remark}
In the case of steady states, the attractor $\mathcal{C}=\{c\}$ may be denoted simply by $c$.
\end{remark}

\begin{example}\label{eg:cycle_structure}
The Boolean network from Example~\ref{eg:wd} %from the previous example, $F(x_1,x_2,x_3)=(x_2 \wedge \neg x_3,x_3,\neg x_1 \wedge x_2)$, 
has 2 steady states (i.e., attractors of length 1) and one cycle of length 2, which can be easily derived from its state space representation (Figure~\ref{fig:dynamics_eg}b).
\end{example}

%\section{Modules and Combinations of Modules}
\section{Simple networks}

In this section, we define a subnetwork of a Boolean network as the restriction of that network to the variables of a subgraph of the wiring graph.  Of importance will be those subnetworks whose wiring graph is strongly connected.  These will be our simple networks, and we will show how multiple networks can be combined into a larger network. To begin, we introduce the idea of restrictions of Boolean functions and Boolean networks.

\subsection{Restrictions of functions and networks}\label{section:restriction}
In order to define a subnetwork, we must have a prescribed definition of restriction.  In general, one might define the restriction of a function $f$ to a subset of its variables in a variety of ways. Each of these choices however may lead to a different subnetwork, with different dynamics. 

%In order to define a subnetwork, we must have a clear definition of restriction.  In general, one might define the restriction of a function $f$ to a subset of its variables in a variety of ways. Each of these choices however may lead to a variety of different dynamics of the restricted subnetwork.  

\begin{example}
Consider the function
    \[f(x_1,x_2,x_3)=(x_1\wedge x_2)\vee(\neg x_1\wedge x_3)=\begin{cases}
        x_2 & \text{if }x_1=1,\\
        x_3 & \text{if }x_1=0.\\
    \end{cases}\]
In order to restrict $f$ to the set $Y=\{x_2,x_3\},$ we have two choices.  We may either set $x_1=1$ in which case $f|_Y=x_2$, or we can set $x_1=0$ yielding a different function $f|_Y=x_3.$ 
%Yet, another way to restrict $f$ to $Y$ can be obtained by setting $x_1=1$ in $x_1\wedge x_2$ and $x_1=0$ in $\neg x_1\wedge x_3$, in this case $f_Y=x_2\vee x_3$.
\end{example}
This example highlights that, in general, the restricted subnetwork depends on the specific choice of restriction and there is no clear choice for this restriction. For this reason, we consider here the class of fanout-free functions, which comes with a natural, unambiguous choice of restriction. 
\begin{definition}\label{Subnetwork}
Given a fanout-free function $h:\FF_2^n\to\FF_2$ and a subset of its variables $Y\subset \{x_1,\ldots,x_n\}$, we define the {\it restriction of $h$ to $Y$} to be the function $h|_Y$ where for all $x_j\notin Y$ we set $x_i$ as follows:
\begin{enumerate} 
\item If $h$ has the form $h=\cdots(x_j\vee\cdots)\cdots\quad$ or $h=\cdots(\neg x_j\wedge\cdots)\cdots\quad,$ then set $x_j=0$;
\item If $h$ has the form $h=\cdots(x_j\wedge\cdots)\cdots\quad$ or $h=\cdots(\neg x_j\vee\cdots)\cdots\quad,$ then set $x_j=1$.
\end{enumerate}
If $h$ depends only on a single variable, i.e., $h=x_j$ or $h=\neg x_j$, and $x_j \not\in Y$, then we set $h|_Y=1$. Note that this arbitrary choice of the constant agrees with the treatment of this case in exceptional case 1 in Theorem~\ref{thm:he}.
\end{definition}

\begin{example}
Consider the fanout-free function
\[h(x_1,x_2,x_3,x_4)=(x_1\wedge \neg x_2)\vee (x_3 \wedge x_4).\] 
If $Y=\{x_1,x_2,x_3\},$ then $h|_Y(x_1,x_2,x_3)=(x_1\wedge \neg x_2)\vee (x_3 \wedge 1)=(x_1\wedge \neg x_2)\vee x_3$. If $Y=\{x_1,x_4\}$, then $h|_Y(x_1,x_4)=(x_1\wedge \neg 0)\vee (1 \wedge x_4)=x_1\vee x_4$. Note that the restriction acts in a way that removal of the variables not in $Y$ does not alter the dependency of $h|_Y$ on all variables in $Y$.
\end{example}

As nested canalizing functions make up a subset of fanout-free functions, they inherit this natural choice of restriction.  Using the rules of Defintion \ref{Subnetwork}, the restriction of a nested canalizing function $f$ to $Y\subset \{x_1,\ldots,x_n\}$ is equivalent to setting all variables not in $Y$ to their non-canalizing input values.  

\begin{example}
Given the nested canalizing function
    \[f(x_1,x_2,x_3)=\begin{cases}
                        b_{1} & x_1 = a_1,\\
b_{2} & x_1 \neq a_1, x_2 = a_2,\\
b_{3} & x_1 \neq a_1, x_2 \neq a_2, x_3 = a_3,\\
1 + b_{3} & x_1 \neq a_1, x_2 \neq a_2, x_3 \neq a_3,\\
                     \end{cases}
                     \]
and the subset $Y=\{x_1,x_3\},$ the restriction of $f$ to $Y$ is given by
    \[f|_Y(x_1,x_3)=f(x_1,a_2+1,x_3)=
    \begin{cases}
        b_{1} & x_1 = a_1,\\
        b_{3} & x_1 \neq a_1, x_3 = a_3,\\
        1 + b_{3} & x_1 \neq a_1, x_3 \neq a_3.
        \end{cases}\]
\end{example}
Now that we have a notion for the restriction of Boolean functions, we can extend that notion to Boolean networks as a whole.
\begin{definition}\label{ex_Subnetwork}
For a Boolean network $F=(f_1,\ldots,f_n):\FF_2^n\to\FF_2^n$ and $Y\subseteq\{x_1,\ldots,x_n\}$, where without loss of generality we consider $Y=\{x_1,\ldots,x_k\}$, we define the \emph{restriction of $F$ to $Y$} by $(F|_Y)_i=(f_i)|_Y$ for $i=1,\ldots,k.$
\end{definition}

\begin{example}
Consider the Boolean network 
\[F(x)=(x_2\wedge x_1, \neg x_1, x_1\vee \neg x_4, (x_1\wedge \neg x_2)\vee (x_3\wedge x_4))\] with wiring diagram in Figure~\ref{fig:modules_eg}a. The restriction of this network to $Y=\{x_1,x_2\}$ is the 2-variable network $F|_Y(x_1,x_2)=(x_2\wedge x_1, \neg x_1)$ with wiring diagram in Figure~\ref{fig:modules_eg}b, while the restriction of $F$ to $Y=\{x_3,x_4\}$ is the 2-variable network $F|_Y(x_3,x_4)=(\neg x_4, x_3\wedge x_4 )$ with wiring diagram in Figure~\ref{fig:modules_eg}c. Note that the wiring diagram of $F|_Y$ is always a subgraph of the wiring diagram of $F$, irrespective of the choice of $Y$.
\end{example}

% \begin{figure}
%     \centering
%     \includegraphics[scale=1.5]{restriction_eg.pdf}
%     \caption{Restriction of Boolean networks. 
%     (a) Wiring diagram of Boolean network $F$. (b) Wiring diagram of Boolean network $F$ restricted to $\{x_3,x_4\}$, which will always be a subgraph of the full wiring diagram.
%     }
%     \label{fig:restriction_eg}
% \end{figure} 

\subsection{Simple Networks}\label{subsec:simple_networks}
Having defined restrictions for the class fanout-free functions, we can now describe the elementary components that will be employed in our decomposition theory in Section~\ref{sec:decomposition}.  We will call these components \emph{simple networks}, in analogy to the decomposition of a finite group into simple groups. These simple networks will be determined solely by the wiring diagram.

\begin{definition}\label{def_strongly_connected}
The wiring diagram of a Boolean network is \emph{strongly connected} if every pair of nodes is connected by a directed path. That is, for each pair of nodes $x_i,x_j$ in the wiring diagram  with $x_i\neq x_j$ there exists a directed path from $x_i$ to $x_j$ (and vice versa). In particular, a one-node wiring diagram is strongly connected by definition.
\end{definition}

\begin{figure}
    \centering
    \includegraphics[scale=1.5]{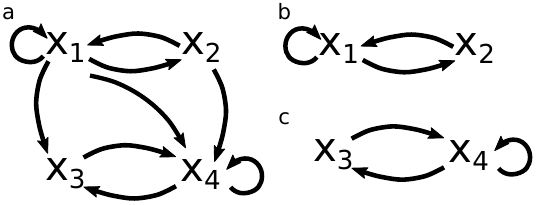}
    \caption{Boolean network decomposition. (a) Wiring diagram of a non-strongly connected Boolean network $F$. (b-c) Wiring diagram of $F$ restricted to (b) $\{x_1,x_2\}$ and (c) $\{x_3,x_4\}$. These strongly connected components are the wiring diagrams of the simple networks of $F$.
    }
    \label{fig:modules_eg}
\end{figure}

\begin{example}
The wiring diagram in Figure~\ref{fig:modules_eg}a is not strongly connected but the wiring diagrams in Figure~\ref{fig:modules_eg}bc are strongly connected.
\end{example}

%For clarity, our definitions and results will be made for networks where each variable is part of a directed cycle. Since such nodes can be removed without changing the dynamical properties \cite{Saadatpour2013reduction}, this will not cause any loss of generality.\claus{@ALAN: Why do we need to make this statement? Can we delete it? If not, we need to explain better what we mean by the first sentence}

\begin{definition}\label{def_simple}
A Boolean network $F:\FF_2^n\to \FF_2^n$ is \emph{simple} if its wiring diagram is strongly connected. If $n=1$ and $F$'s wiring diagram has no edge (i.e., no self-regulation), we call the simple network $F$ \emph{trivial}.
\end{definition}

\begin{remark}
For an arbitrary Boolean network $F$, its wiring diagram will either be strongly connected or it will be made up of a collection of strongly connected components where connections between each component move in only one direction.  

Let $F$ be a Boolean network and let $W_1,\ldots,W_m$ be the strongly connected components of its wiring diagram, with $X_i$ denoting the set of variables in strongly connected component $W_i$. Then, the \emph{simple subnetworks} of $F$ are defined by $F|_{X_1},\ldots,F|_{X_n}$.
\end{remark}

\begin{definition}\label{def:acyclic}
Let $W_1,\ldots, W_m$ be the strongly connected components of the wiring diagram of a Boolean network $F$. By setting $W_i \rightarrow W_j$ if there exists at least one edge from a vertex in $W_i$ to a vertex in $W_j$, we obtain a (directed) acyclic graph 
$$Q = \{(i,j) | W_i \rightarrow W_j \},$$ which describes the connections between the strongly connected components of $F$. %Note that this acyclic graph induces a partial order on the modules    (see \cite{jarrah2010dynamics}).

%We define a partial order by setting $W_i \preccurlyeq W_j$ if there exists at least one edge from a vertex in $W_i$ to a vertex in $W_j$. The set of modules together with $\preccurlyeq$ forms a partially ordered set%, which can be expressed as a lower triangular matrix $Q$ with
% $$q_{ij} = \begin{cases}
%     1 & \ \text{if }W_i \preccurlyeq W_j,\\
%     0 & \ \text{otherwise.}
% \end{cases}$$
\end{definition}

%\claus{add example of trivial modular network?}

\begin{example}
For the Boolean network $F$ from Example~\ref{eg:wd}, the wiring diagram
%For the previous example, the wiring diagram (Figure~\ref{fig:modules_eg}a) of the Boolean network \[F(x)=(x_2\wedge x_1, \neg x_1, x_1\vee \neg x_4, (x_1\wedge \neg x_2)\vee (x_3\wedge x_4))\]
has two strongly connected components $W_1$ and $W_2$ with variables $X_1=\{x_1,x_2\}$ and $X_2=\{x_3,x_4\}$ (Figure~\ref{fig:modules_eg}bc), connected according to the directed acyclic graph $Q = \{(1,2)\}$. The two simple networks of $F$ are given by the restriction of $F$ to $X_1$ and $X_2$, that is, $F|_{X_1}(x_1,x_2) = (x_2\wedge x_1, \neg x_1)$ and $F|_{X_2}(x_3,x_4) = (\neg x_4, x_3\wedge x_4)$. 
%Note that since the ``upstream component'' (with variables $x_1$ and $x_2$) does not receive feedback from the other component, the corresponding simple network (i.e., the restriction of $F$ to $X_1$) is simply the projection of $F$ onto the variables $(x_1,x_2)$.  
Note that the simple network $F|_{X_1}$, i.e., the restriction of $F$ to $X_1$, is simply the projection of $F$ onto the variables $X_1$ because $W_1$ does not receive feedback from the other component (i.e., because $(2,1)\not\in Q$).
\end{example}

%So far we have described how to formally define and describe the modules of any Boolean network governed by fanout-free functions. We consider fanout-free functions because they allow for a natural way to define the restriction to a set of variables. However, as long as there is an unambiguous way to define such a restriction, all following results still apply. Linear functions constitute another family of Boolean functions where there is a natural definition of restriction. For example, the restriction of $f(x_1,x_2,x_3,x_4)=x_1+x_2+x_3+x_4$ to $Y=\{x_1,x_4\}$ is simply $f|_Y(x_1,x_4)=x_1+x_4$.

%We now define a semi-direct product that combines networks governed by functions with an unambiguously defined restriction.

So far we have described how to formally define and describe restrictions and simple networks. We now define a semi-direct product that combines two networks. 

\begin{definition}
Consider two Boolean networks, $F=(f_1,\ldots,f_k):\FF_2^k\to\FF_2^k$ with variables $\{x_1,\ldots,x_k\}$ and $G=(g_1,\ldots,g_r):\FF_2^r\to\FF_2^r$ with variables $\{y_1,\ldots,y_r\}$. 
%with $W_1:=\{x_1,\ldots,x_k\}$ and $W_2:=\{y_{1},\ldots,y_{r}\}.$.
 Let $P=\{(y_i,p_i)\}_{i\in\Lambda}$ where $\Lambda\subset\{1,\ldots,r\}$, $p_i: \FF_2^{k+r}\to\FF_2$, and $(p_i)|_{X_2}=g_i.$ Then we define a new Boolean network $H=(h_1,\ldots,h_{k+r}):\FF_2^{k+r}\to \FF_2^{k+r}$ by 
$$h_i(x,y)=f_i(x) \text{ for } 1\leq i\leq k$$
and 
$$
h_{k+i}(x,y)=\begin{cases}
p_i(x,y) & \text{if }i\in\Lambda,\\
g_i(y) & \text{if }i\not\in\Lambda
\end{cases}
$$
for the other subindices.

We denote $H$ as $H:=F\rtimes_P G$ and refer to this as a coupling of $F$ and $G$ by (the coupling scheme) $P$ or as the semi-direct product of $F$ and $G$ via $P$.  If $P$ is clear from the context (e.g., if the full network is given explicitly as in the previous example), we can omit $P$ and just use the notation $F\rtimes G$.
\end{definition}

\begin{remark}\label{remark:any_fun_works}
Note that, by definition, $F\rtimes_P G$ restricted to $X_1$ and $X_2$ results in $F$ and $G$; that is, 
$(F\rtimes_P G)|_{X_1} = F$ and $(F\rtimes_P G)|_{X_2} = G$. This product is defined as long as there is a prescribed way to define a restriction.
We have shown that fanout-free functions allow for a natural way to define the restriction to a set of variables (Definition~\ref{Subnetwork}). 
%However, as long as there is an unambiguous way to define such a restriction, all following results still apply. 
Linear functions constitute another family of Boolean functions where there is a natural definition of restriction. For example, the restriction of $f(x_1,x_2,x_3,x_4)=x_1+x_2+x_3+x_4$ to $Y=\{x_1,x_4\}$ is simply $f|_Y(x_1,x_4)=x_1+x_4$.
\end{remark}

\begin{figure}
    \centering
    \includegraphics[scale=1.5]{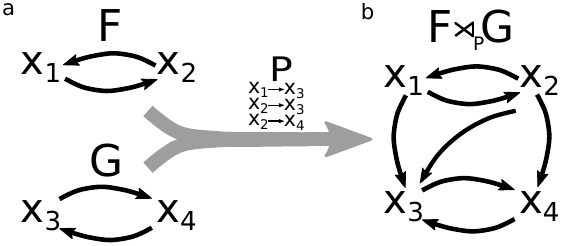}
    \caption{semi-direct product of Boolean networks. (a) Independent Boolean networks $F$ and $G$ and a coupling scheme $P$, which determines which edges to add between the wiring diagrams of $F$ and $G$. (b) The resulting network is a semi-direct product of $F$ and $G$ by (the coupling scheme) $P$. 
    }
    \label{fig:product_eg}
\end{figure}

\begin{example}
Consider the Boolean networks %\[F(x_1,x_2)=(x_2,x_1)\quad\text{and}\quad G(x_3,x_4)=(x_4,x_3)\] (or equivalently as in the above definition, $G(y_1,y_2)=(y_2,y_1)$).
\[F(x_1,x_2)=(x_2,x_1)\quad\text{and}\quad G(y_1,y_2)=(y_2,y_1)\]
If we use the coupling scheme $P=\{(y_1,x_1\vee(x_2\wedge y_2)),(y_2,\neg x_2 \wedge y_1 )\}$, we obtain the nested canalizing network $F \rtimes_P G: \FF_2^4\to \FF_2^4$,
$$(F\rtimes_P G)(x_1,x_2,y_1,y_2)=(x_2,x_1,x_1\vee(x_2\wedge y_2),\neg x_2 \wedge y_1).$$ 
At the wiring diagram level, this product can be seen as the union of the two wiring diagrams and some added edges determined by the coupling scheme $P$ (Figure~\ref{fig:product_eg}).

%For instance, $(x_3,x_1\vee(x_2\wedge x_4))$ indicates that we need to add (positive) edges from $x_1$ and $x_2$ to $x_3$ (the edge from $x_4$ is already present); $(x_4,\neg x_2 \wedge x_3 )$ indicates that we need to add a (negative) edge from $x_2$ to $x_4$.
\end{example}

\begin{example}\label{example_linear}
If $F(x_1,x_2)=(x_2,x_1)$ and $G(y_1,y_2)=(y_2,y_1)$, as before, and we use the coupling scheme $P=\{(y_1,x_1+x_2+y_2),(y_2, x_2 + y_1 )\}$, then we obtain the linear network
$$(F\rtimes_P G)(x_1,x_2,y_1,y_2)=(x_2,x_1,x_1+x_2+y_2,x_2 + y_1).$$
At the wiring diagram level, this product also looks like Figure~\ref{fig:product_eg}.
\end{example}

%\begin{remark}
%This definition directly implies the following: A Boolean network is modular if and only if it is not strongly connected.
%\end{remark}

%\subsection{Twisted and non-autonomous networks}

%\section{Existence and Uniqueness}% - ALAN, maybe we should include more from above in this section

\subsection{Cross product of networks} 
When $P= \emptyset$, %, $\Lambda$ will be empty and %the second case in the definition of $H$ will not occur; then 
$F_1 \rtimes_{\{\}} F_2=F_1\times F_2$ is simply a Cartesian product of networks. In this case the dynamics of $F =F_1\times F_2 $ can be determined directly from the dynamics of $F_1$ and $F_2$.

Consider a network $F:\FF_2^n\rightarrow\FF_2^n$ such that the variables of $F$ can be partitioned so that $F(x,y)=(F_1(x),F_2(y))$ for some Boolean networks $F_1:\FF_2^{n_1}\rightarrow \FF_2^{n_1}$ and $F_2:\FF_2^{n_2}\rightarrow \FF_2^{n_2}$ with $n=n_1+n_2$. In this case, we can easily decompose $F$ as a Cartesian product, $F=F_1\times F_2$. Also, the wiring diagram of $F$ is the disconnected union of the wiring diagrams of $F_1$ and $F_2$.

The dynamics of $F$ consists of coordinate pairs $(x,y)$ such that 
\[x(t+1)=F_1(x(t)), \quad y(t+1)=F_2(y(t)).\]
Furthermore, the trajectory $(x(t),y(t))_{t=0}^\infty$ of $F$ is periodic if and only if  $(x(t))_{t=0}^\infty$ and $(y(t))_{t=0}^\infty$ are periodic trajectories of $F_1$ and $F_2$, respectively. This provides a way to decompose the dynamics of $F$ as the product of the dynamics of $F_1$ and $F_2$.

\begin{lemma}\label{lem:lcm}
If trajectories $(x(t))_{t=0}^\infty$ and $(y(t))_{t=0}^\infty$ have periods $m$ and $l$, respectively, then the trajectory $(x(t),y(t))_{t=0}^\infty$ has period $lcm(m,l)$. 
\end{lemma}

\begin{proposition}
The set of periodic points of $F$ is the Cartesian product of the set of periodic points of $F_1$ and periodic points of $F_2$.
\end{proposition}
%\begin{proof}
%It is enough to note that $(x,y)$ is a periodic point of $F$ if and only if $x$ is a periodic point of $F_1$ and $y$ is a periodic point of $F_2$. Furthermore, if $p_1$ is the period of $x$ and $p_2$ is the period of $y$, then $lcm(p_1,p_2)$ is the period of $(x,y)$.
%\end{proof}

\begin{figure}
    \centering
    \includegraphics[scale=0.3]{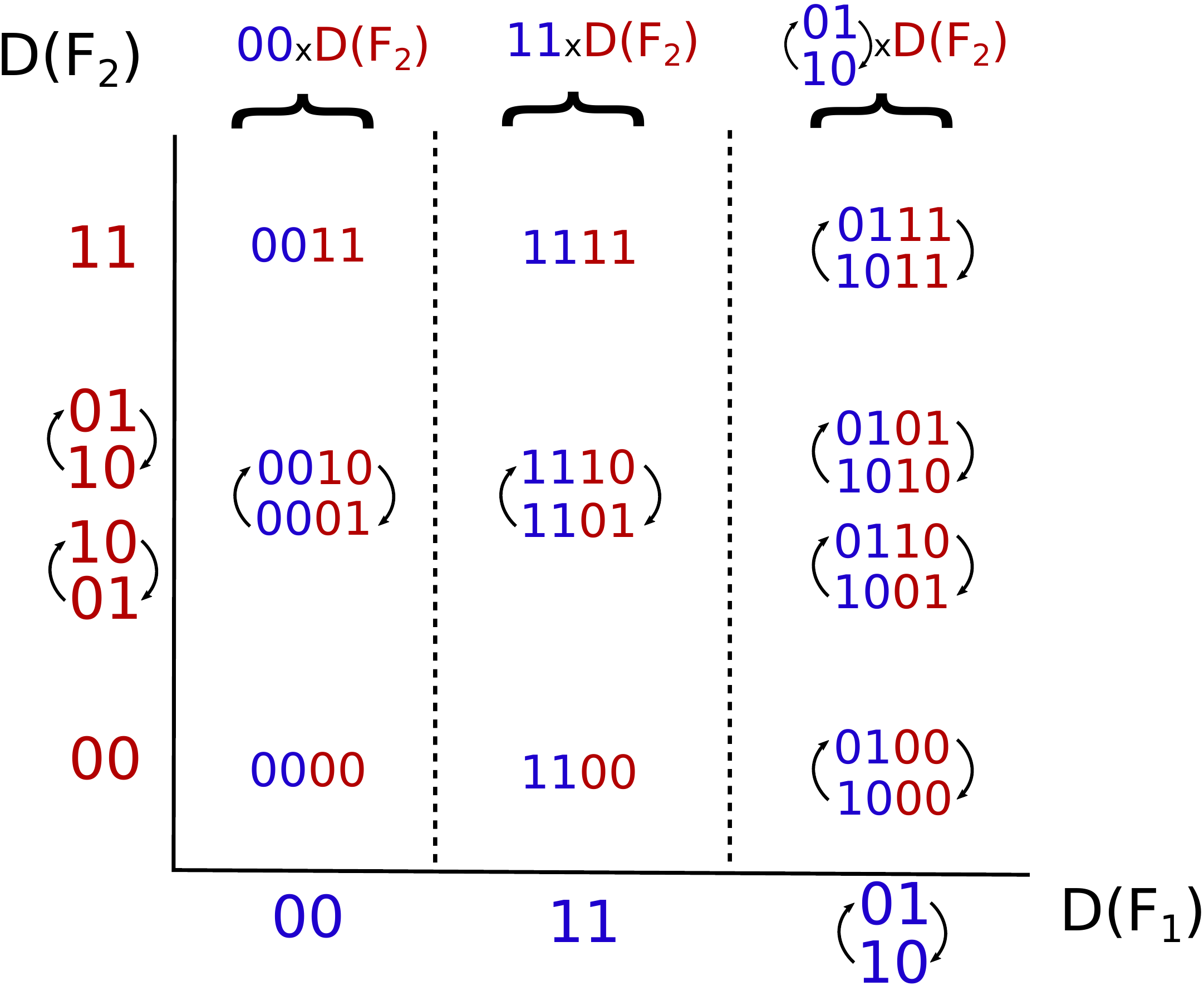}
    \includegraphics[scale=0.3]{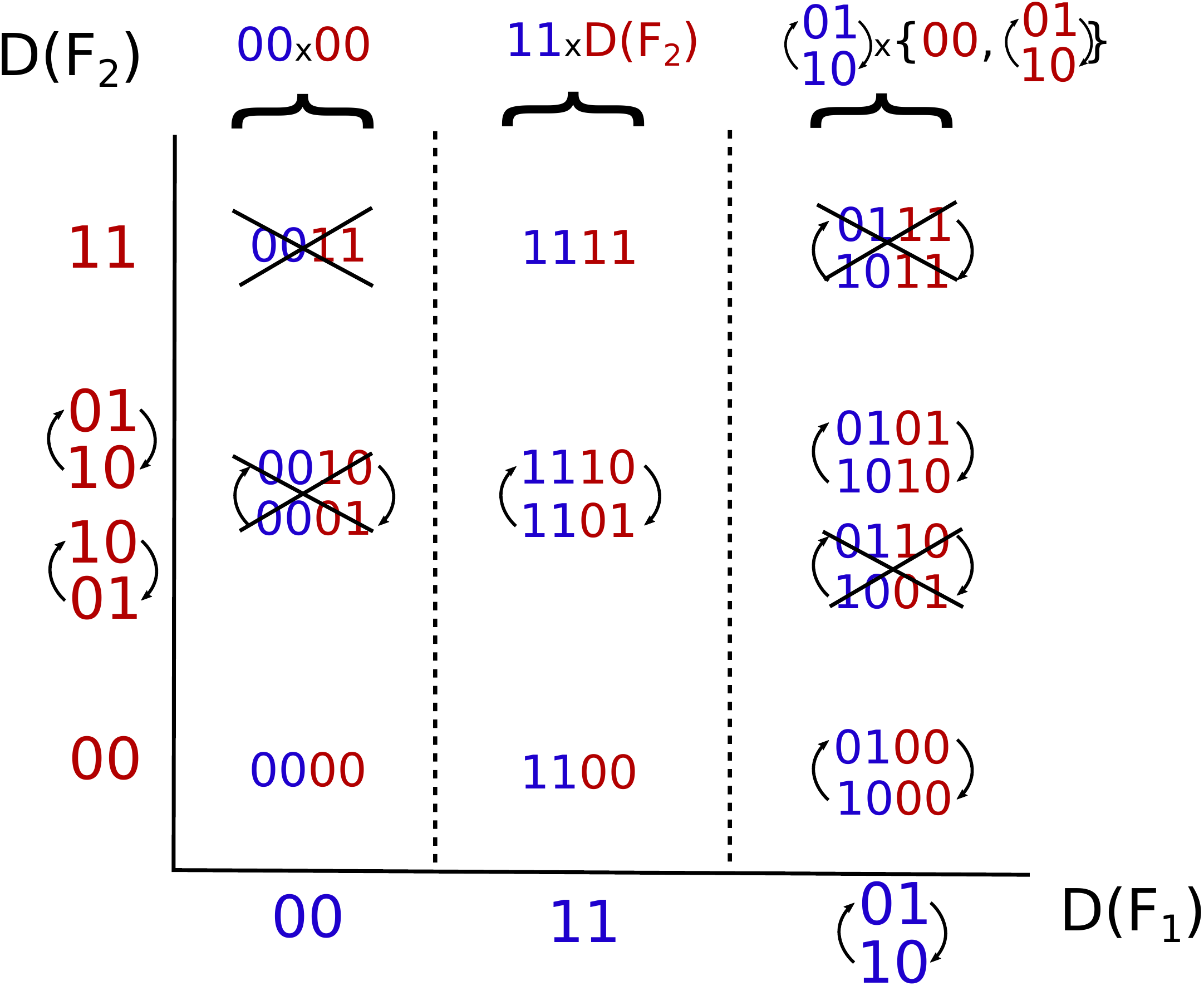}
    \caption{Dynamics of a Cartesian product and a semi-direct product.  
    \textbf{Left:} The dynamics of a Cartesian product $F=F_1\times F_2$ can be seen as a Cartesian product of $\mathcal{D}(F_1)$ and $\mathcal{D}(F_2)$.
    \textbf{Right:} The coupling of networks does not behave as a Cartesian product and the dynamics can vary due to this coupling. In this example, the crossed-out attractors indicate which attractors from the Cartesian product are lost when using a semi-direct product.
    }
    \label{fig:cart_prod}
    \end{figure}

\begin{example}
Consider $F:\FF_2^4\rightarrow\FF_2^4$ given by $F(x_1,x_2,x_3,x_4)=(x_2,x_1,x_4,x_3)$. $F$ can be seen as $F=F_1\times F_2$, where $F_1(x_1,x_2)=(x_2,x_1)$ and $F_2(x_3,x_4)=(x_4,x_3)$. The dynamics of $F_1$, $F_2$, and $F$ are shown in the left panel of Figure~\ref{fig:cart_prod}.
The dynamics of $F_1$ is $\mathcal{D}(F_1)=\{00,11,\{01,10\}\}$ and the dynamics of $F_2$ is $\mathcal{D}(F_2)=\{00,11,\{01,10\}\}$ (note that we are denoting steady states $\mathcal{C}=\{c\}$ simply by $c$). By concatenating the attractors of $F_1$ and $F_2$, we obtain the attractors of $F$. Note that we have two ways of concatenating $\{01,10\}$ (attractor of $F_1$) and $\{01,10\}$ (attractor of $F_2$) to obtain attractors of $F$. If we start at 0101 we obtain the attractor $\{0101,1010\}$, and if we start at $0110$ we obtain the attractor $\{0110,1001\}$. This difference is captured by explicitly including both, $\{01,10\}$ and $\{10,01\}$, as periodic trajectories of $F_2$ (vertical axis). 
In this case, we can see that $\mathcal{D}(F)$ is the union of $00\times \mathcal{D}(F_2)$, $11\times \mathcal{D}(F_2)$, and $\{01,10\}\times \mathcal{D}(F_2)$.  
\end{example}

In the general case, $F$ will not be a Cartesian product, but can be seen as the coupling of two or more networks.

\begin{example}\label{example:small_coupled}
Consider $F:\FF_2^4\rightarrow\FF_2^4$ given by $F(x_1,x_2,x_3,x_4)=(x_2,x_1,x_2 x_4,x_3)$. $F$ can be seen as $F=F_1\rtimes_P F_2$, where $F_1(x_1,x_2)=(x_2,x_1)$ and $F_2(x_3,x_4)=(x_4,x_3)$, as before, but $P=\{(x_3,x_2 x_4)\}$ is no longer empty. The dynamics of $F_1,F_2$ and $F$ are shown in right panel of Figure~\ref{fig:cart_prod}.
The dynamics of $F_1$ is $\mathcal{D}(F_1)=\{00,11,\{01,10\}\}$  and the dynamics of $F_2$ is $\mathcal{D}(F_2)=\{00,11,\{01,10\}\}$. Since $F$ is not the Cartesian product of $F_1$ and $F_2$, we cannot simply take a Cartesian product of dynamics.
Due to the coupling between $F_1$ and $F_2$, not every combination of attractors of $F_1$ and $F_2$ will result in an attractor of $F$. The products that do not result in an attractor of $F$ are crossed out in the figure. Note that $\{01,10\}\in \mathcal{D}(F_1)$ and $\{01,10\}\in \mathcal{D}(F_2)$ \textbf{do} give rise to an attractor of $F$, but $\{01,10\}\in \mathcal{D}(F_1)$ and $\{10,01\}\in \mathcal{D}(F_2)$ \textbf{do not}. 
In this case, we can see that $\mathcal{D}(F)$ is the union of $00\times 00$, $11\times \mathcal{D}(F_2)$, and $\{01,10\}\times \{00,\{01,10\}\}$. Here the dynamics of the coupled network is actually a subset of the dynamics of the Cartesian product. In general, this is not the case but depends on the particular coupling between the networks as will be seen in Section~\ref{sec:dec_dyn}.
\end{example}

%\begin{figure}
 %   \centering
    %\includegraphics[scale=0.3]{dynamics_cartprod.pdf}
 %   \includegraphics[scale=0.3]{dynamics_cartprod_expanded.pdf}
 %   \includegraphics[scale=0.3]{dynamics_cartprod_coupled.pdf}
  %  \caption{Dynamics of a Cartesian product and ``non-Cartesian product''.  
   % \textbf{Left:} The dynamics of a Cartesian product $F=F_1\times F_2$ can be seen as a Cartesian product of $\mathcal{D}(F_1)$ and $\mathcal{D}(F_2)$.
    %
    %\textbf{Right:} The coupling of networks does not behave as a Cartesian product and the dynamics can vary due to this coupling. In this example, some attractors from the Cartesian product are lost.
    %
    %}
    %\label{fig:cart_prod}
    %\end{figure} 

\section{Decomposing networks and network dynamics}\label{sec:decomposition}

\subsection{Decomposing a network}\label{sec:dec_mod}
Example \ref{example:small_coupled} shows that the dynamics of a network can in some ways be seen as a combination of the dynamics of each of its simple networks. Our decomposition theorems in this section formalize this combination. We will show that every network admits a decomposition into a semi-direct product of simple networks, and we formalize how the dynamics of a network can be constructed from this decomposition. In what follows, we will assume that all Boolean networks are fanout-fee but all results still apply in more generality for networks governed by a class of functions with a prescribed restriction (e.g., linear functions, see Remark~\ref{remark:any_fun_works}).

\begin{theorem}\label{thm:modular_exist}
If a Boolean network $F$ is not simple, then there exist $F_1,F_2,P$ such that $F=F_1\rtimes_P F_2$. Furthermore, we can find a decomposition such that $F_1$ is simple.
\end{theorem}
\begin{proof}

%Consider $F=(f_1,\ldots,f_n)$ with variables $x_1,\ldots,x_n$. For a Boolean network to be not strongly connected, it means that there exists at least one node, say $x_0$, and a node, say $x_j$ such that in the wiring diagram of $F$ there does not exist a path from $x_j$ to $x_0$. Let $X_2$ denote the set of all such nodes $X_2=\{x_{j_1},x_{j_2},\ldots,x_{j_m}\}$ for which there are no paths from $x_{j_i}$ to $x_0,$ and let $X_1=X/X_2$ denote the complement set of nodes to $X_2$.
Let $F=(f_1,\ldots,f_n)$ be a Boolean network with variables $X = \{x_1,\ldots,x_n\}$ which is not simple. Then the wiring diagram of $F$ is not strongly connected meaning there exists at least one node $y$ and one node $x_j\neq y$ such that in the wiring diagram of $F$ there exists no path from $x_j$ to $y$. Let $X_2$ denote the set of all such nodes $X_2=\{x_{j_1},x_{j_2},\ldots,x_{j_m}\}$, for which there exists no paths to $y$. Further, let $X_1=X\backslash X_2$ denote the complement set of nodes to $X_2$. Note that for every $x_i\in X_1$, there exists a path from $x_i$ to $y$ but no paths originating from $X_2$ to $x_i$.

%We note here that for any $x_i\in X_1$, there is a path from $x_i$ to $x_0$ and thus no paths originating from any node in $X_2$ and ending at $x_i$.

Define $\Lambda$ to be the subset of indices  $\Lambda\subset\{j_1,\ldots,j_m\}$ such that $f_{j_i}$ depends on at least one $x_i\in X_1$. Note that $\Lambda$ may be empty. Let $F_1=F|_{X_1}$, $F_2=F|_{X_2}$ and $P=\{(x_i,f_i)\}_{i\in \Lambda}$.  We claim that $F=F_1\rtimes_P F_2$.

If $\Lambda = \emptyset$, then the sets $X_1$ and $X_2$ represent two groups of nodes, which are disconnected in the wiring diagram.  Hence the network $F$ is a Cartesian product of $F_1$ and $F_2$.  It follows that $F=F_1\rtimes_P F_2$ where $P=\emptyset.$

Now suppose $\Lambda\neq\emptyset.$ We distinguish three cases to show that $F=F_1\rtimes_P F_2$.
\begin{enumerate}
    \item For any $x_i\in X_1,$ the corresponding update function $(F_1)_i$ does not depend on $X_2$ by construction, as there are no paths from $X_2$ to $x_i$. Thus by definition, $(F_1)_i:=(F|_{X_1})_i=f_i|_{X_1}=f_i.$
    \item For any $x_i\in X_2$ with $i\in\Lambda$, we have $(F_2)_i:=(F|_{X_2})_i=f_i|_{X_2}$ as defined by $P$. 
    \item For any $x_i\in X_2$ with $i\not\in\Lambda$, by definition of $\Lambda$, $f_i$ does not depend on $X_1$. Thus, $(F_2)_i:=(F|_{X_2})_i=f_i|_{X_2}=f_i$.  
\end{enumerate}
Hence it follows that $F=F_1\rtimes_P F_2.$

%If $F_1$ and $F_2$ are not strongly connected, for each $F_i$, we can repeat the above procedure until they are strongly connected. 
Note that in the above proof we can choose the node $y$ such that it belongs to a strongly connected component that receives no edge from any other strongly connected component. $X_1$ will contain the nodes of this strongly connected component and hence $F_1$ will be simple. 
%If the node $y$  is chosen so that it belongs to a strongly connected component that feeds no edge to any other strongly connected component, $X_2$ will be the nodes of a strongly connected component and hence $F_2$ will be strongly connected.
%If the node $y$ is chosen so that it belongs to a strongly connected component that receives no edge from any other strongly connected component, $X_1$ will be the nodes of a strongly connected component and hence $F_1$ will be strongly connected. If the node $y$  is chosen so that it belongs to a strongly connected component that feeds no edge to any other strongly connected component, $X_2$ will be the nodes of a strongly connected component and hence $F_2$ will be strongly connected.
\end{proof}

\begin{definition}\label{rem_modular_network}
A Boolean network $F$ is \emph{decomposable} if there exist networks $F_1$ and $F_2$ and a coupling scheme $P$ such that $F=F_1\rtimes_P F_2$. 
\end{definition}

Recall that each Boolean network has an associated directed acyclic graph $Q$, determined by the connections between strongly connected components of the wiring diagram. This acyclic graph can be extended to a partial order on the set of strongly connected components encoding which nodes are ``upstream" or ``downstream" from a particular module.

% By induction on the downstream component $F_2$ in Theorem~\ref{thm:modular_exist}, we obtain the following result that shows that any fanout-free Boolean network is either simple or decomposable into a series of semi-direct product of simple Boolean networks.

Our first main result shows that any Boolean network is either simple or decomposable into a unique series of semi-direct product of simple Boolean networks.

\begin{theorem}\label{thm:unique_decomp}
If a Boolean network $F$ is not simple, then there exist unique  simple networks $F_1,\ldots,F_m$ such that 
\[F=F_1\rtimes_{P_1}(F_2\rtimes_{P_2}(\cdots\rtimes_{P_{m-1}} F_m)),\]
where this representation is unique up to a reordering, which respects the partial order of $Q$, and the collection of pairings $P_1,\ldots,P_{m-1}$ depends on the particular choice of ordering.
\end{theorem}
\begin{proof}
By induction on the downstream component $F_2$ in Theorem~\ref{thm:modular_exist} we obtain the %following 
result.% that shows that any fanout-free Boolean network is either simple or decomposable into a series of semi-direct product of simple Boolean networks.
\end{proof}
 
\begin{remark}
The parentheses of the representation of $F$ as a decomposition into simple networks may be rearranged in any associative manner.  However each choice determines a new collection of pairings $P'_1,\ldots,P'_{m-1}$.
\end{remark}

The main takeaway is that every Boolean network is either a simple network or it admits a decomposition into a semi-direct product of simple networks.

%\begin{corollary}
%If $F$ is not strongly connected, then there exists unique (up to reordering) strongly connected networks $F_1,\ldots,F_n$ and pairings $P_1,\ldots,P_{n-1}$ such that \[F=((F_1\rtimes_{P_1}F_2)\rtimes_{P_2}\cdots)\rtimes_{P_{n-1}} F_n.\]
%\end{corollary}

\subsection{Non-autonomous networks}
In order to study the coupling of networks, we need to understand how a trajectory, which will describe an attractor of an ``upstream" network, influences the dynamics of a ``downstream" network. This can be studied using the concept of non-autonomous networks.
\begin{definition}
A \emph{non-autonomous Boolean network} is defined by $$y(t+1)=H(g(t),y(t)),$$
where $H:\FF_2^{m+n}\rightarrow \FF_2^n$ and $\left(g(t)\right)_{t=0}^\infty$ is a sequence with elements in $\FF_2^m$. We call this type of network non-autonomous because its dynamics will depend on $g(t)$. We use $H^g$ to denote this non-autonomous network.

A state $c\in\FF_2^n$ is a \emph{steady state} of $H^g$ if $H(g(t),c)=c$ for all $t$. Similarly, an ordered set with $r$ elements, $\mathcal{C}=\{c_1,\ldots,c_r\}$ is an \emph{attractor of length $r$} of $H^g$ if $c_2=H(g(1),c_1)$, $c_3=H(g(2),c_2),\ldots, 
c_r=H(g(r-1),c_{r-1})$,
$c_1=H(g(r),c_r)$, $c_2=H(g(r+1),c_1),\ldots$. Note that in general $g(t)$ is not necessarily of period $r$ and may even not be periodic.
\end{definition}

If $H(g(t),y)=G(y)$ for some network $G$ (that is, it does not depend on $x$) for all $t$, then $y(t+1)=H(g(t),y(t))=G(y(t))$  and this definition of attractors coincides with the classical definition of attractors for (autonomous) Boolean networks (Definition~\ref{def_attractors}).

%With this notation, if $g$ ``covers'' an attractor $\mathcal{C}$ (ie. $\mathcal{C}=\{g(0),g(1),\ldots\}$), then $\bar F_2^\mathcal{C} $ ``is'' the non-autonomous system given by $F_2(g(t),x_2)$. [ Need to add a note that the attractors do not depend on where $g$ begins.]

\begin{example}
Consider the non-autonomous network defined by \[H(x_1,x_2,y_1,y_2)=(x_2 y_2, y_1)\] 
and the two-periodic sequence $\left(g(t)\right)_{t=0}^\infty=(01,10,01,10,\ldots)$. If the initial point is $y(0)=(y_1^*,y_2^*)$, then the dynamics of $H^g$ can be computed as follows: 
\begin{align*}
y(1)&=H(g(0),y(0))=H(0,1,y_1^*,y_2^*)
=(y_2^*,y_1^*),\\
y(2)&=H(g(1),y(1))
=H(1,0,y_2^*,y_1^*)=(0,y_2^*),\\
y(3)&=H(g(2),y(2))=H(0,1,0,y_2^*)=(y_2^*,0).
\end{align*} 
Thus for $t\geq 1$, $y(2t)=(0,y_2^*)$ and $y(2t+1)=(y_2^*,0)$. It follows that the attractors of $H^g$ are given by $\{00\}$ (one steady state) and $\{01,10\}$ (one cycle of length $2$). Note that $\{10,01\}$ is not an attractor because $(10,01,10,01,...)$ is not a trajectory for this non-autonomous network. This is a subtle situation that can be sometimes missed when using the compact representation of an attractor as a set, but is explicitly handled by the trajectory representation.
\end{example}

\begin{example}
Consider the non-autonomous network defined by $H(x_1,x_2,y_1,y_2)=(x_2 y_2, y_1)$, as in the previous example, and the one-periodic sequence $\left(g(t)\right)_{t=0}^\infty=(00,00,\ldots)$. If the initial point is $y(0)=(y_1^*,y_2^*)$, then the dynamics of $H^g$ can be computed as follows: \begin{align*}
y(1)&=H(g(0),y(0))=H(0,0,y_1^*,y_2^*)
=(0,y_1^*),\\
y(2)&=H(g(1),y(1))
=H(0,0,y_2^*,y_1^*)=(0,0).
\end{align*}
Then,
$y(t)=(0,0)$ for $t\geq 2$, and the only attractor of $H^g$ is the steady state $\{00\}$.
\end{example}

\subsection{Decomposing dynamics}\label{sec:dec_dyn}

For a decomposable network $F:=F_1\rtimes_P F_2$, we introduce the following notation for attractors. First, note that $F$ has the form $F(x,y)=(F_1(x),\overline{F}_2(x,y))$ where $\overline{F_2}$ is a non-autonomous network. Let $\cC_1=\{r_1,\ldots,r_m\}\in\cD(F_1)$ and $\cC_2=\{s_1,\ldots, s_n\} \in \cD(F_2^{\cC_1}) := \cD(\overline{F}_2^{\cC_1})$ be attractors. By Lemma~\ref{lem:lcm}, we know that the sequence $((r_t,s_t))_{t=0}^\infty$ has period $l=lcm(m,n)$, so we define the sum or concatenation of these attractors to be
%\[\cC_1\oplus\cC_2=((r_1,s_1),(r_2,s_2),\ldots,(r_m,s_m),(r_{1},s_{m+1}),\ldots,(r_m,s_{2m}),(r_m,s_{2m+1}),\ldots,(r_m,s_n)).\]
\[\cC_1\oplus\cC_2=((r_1,s_1),(r_2,s_2),\ldots,(r_{l-1},s_{l-1})).\]
Note that the sum of attractors is not  a  Cartesian product, $\mathcal{C}_1\times \mathcal{C}_2=\{(r_i,s_j)|\text{ for all } i,j \}$.

Similarly, for an attractor $\mathcal{C}_1$ and a collection of attractors $D$ we define 
\[
\mathcal{C}_1\oplus D=\{\mathcal{C}_1 \oplus \mathcal{C}_2 | \mathcal{C}_2\in D\}.
\]
To encode all $\mathcal{D}(F^{\mathcal{C}_1}_2)$'s where $\mathcal{C}_1\in \mathcal{D}(F_1)$, we use the notation 
\[\mathcal{D}(F_2 ^{F_1}) := \left( \mathcal{D}(F^{\mathcal{C}_1}_2)\ \right )_{\mathcal{C}_1\in\mathcal{D}(F_1)}\]
and refer to this as the  \emph{dynamics of $F_2$ induced by $F_1$}.

% \begin{theorem}\label{thm:decomp}
% Let $F=F_1\rtimes_P F_2$ be a decomposable network.  Then \[\mathcal{D}(F)
% = \bigsqcup_{\mathcal{C}_1\in\mathcal{D}(F_1)}\mathcal{C}_1\oplus\mathcal{D}(F^{\mathcal{C}_1}_2)
% =:\mathcal{D}(F_1)\rtimes \mathcal{D}(F_2 ^{F_1}) 
% .\]
% \claus{we use $oplus$ here but $\times$ in the two figures}
% That is, 
% \[\mathcal{D}(F_1\rtimes_P F_2)=\mathcal{D}(F_1)\rtimes \mathcal{D}(F_2 ^{F_1})\]
% which can be seen as a distributive property for dynamics of coupled networks.

% \end{theorem}

Our second main result shows that the dynamics of a semi-direct product can be seen as a type of semi-direct product of the dynamics of the simple networks.

\begin{theorem}\label{thm:decomp}
Let $F=F_1\rtimes_P F_2$ be a decomposable network.  Then \[\mathcal{D}(F)
= \bigsqcup_{\mathcal{C}_1\in\mathcal{D}(F_1)}\mathcal{C}_1\oplus\mathcal{D}(F^{\mathcal{C}_1}_2).\]
\end{theorem}
\begin{proof}
Let $\cC=\{c_1,\ldots,c_n\}\in\cD(F)$. Let $\cC_1=\text{pr}_1(\cC)=\{\text{pr}_1(c_1),\ldots,\text{pr}_1(c_n)\}=:\{c_1^1,\ldots,c_n^1\}$ and similarly $\cC_2=\text{pr}_2(\cC)=:\{c_1^2,\ldots,c_n^2\}.$ Since $F_1$ by definition does not depend in any way on $X_2$, then  $F_1(\text{pr}_1(x))=\text{pr}_1(F(x)).$ Thus for any $c_j^1$,
\[F_1(c_j^1)=F_1(\text{pr}_1(c_j))=\text{pr}_1(F(c_j))=\text{pr}_1(c_{j+1})=c_{j+1}^1.\]
Iterating this, we find that in general $F_1^k(c_j^1)=c_{j+k}^1$ from which it follows that $\cC_1\in\cD(F_1)$.\\

Next, we consider the non-autonomous network $\overline{F}_2^{\mathcal{C}_1}$ defined as above where $y(t+1)=\text{pr}_2F(g(t),y(t)),$and  $g(t)=c_t^1$.  If $y(1)=c_1^2$, then
    \[y(2)=\pb F(g(1),c_1^2)=\pb F(c_1^1,c_1^2)=\text{pr}_2F(c_1)=\text{pr}_2(c_2)=c_2^2\] 
and in general
    \[y(k+1)=\pb F(g(k),c_k^2)=\pb F(c^1_k,c^2_k)=\pb c_{k+1}=c_{k+1}^2\]
Hence $y(n+1)=\pb F(c_n)=\pb c_1=y(1)$ and thus $\cC_2\in\cD(F_2^{\mathcal{C}_1}).$ From this we have that $\cC=\cC_1\oplus\cC_2\in\cC_1\oplus\cD(F_2^{\mathcal{C}_1})$ and thus \[\cD(F)\subset\bigsqcup_{\cC_1\in\cD(F)}\cC_1\oplus\cD(F_2^{\mathcal{C}_1}).\]

Conversely, let $\cC_1\in\cD(F_1)$ and $\cC_2\in\cD(F_2^{\mathcal{C}_1})$.  We want to show that $\cC_1\oplus\cC_2\in\cD(F)$. Let $g(t)=c^1_t,\,y(1)=c_1^2$, and $y(t+1)=\pb F(g(t),y(t)).$ Since $\cC_2\in\cD(F_2^{\mathcal{C}_1}),$ then $y(t+1)=c_{t+1}^2$ by definition.  Let $N=|\cC_2|$.  Then
    \[F(c_k^1,c_k^2)=(\pa F(c_k^1,c_k^2),\pb F(g(k),y(k))=(F_1(c_k^1),\overline{F}_2^{\mathcal{C}_1}(c_k^1,y(k+1)))=(c_{k+1}^1,c_{k+1}^2).\]
Thus $F^N(c_1^1,c_1^2)=F(c_N^1,c_N^2)=(c_1^1,c_1^2)$ and hence $\cC_1\oplus\cC_2\in\cD(F).$ It follows that \[\bigsqcup_{\cC_1\in\cD(F_1)}\cC_1\oplus\cD(F_2^{\mathcal{C}_1})\subset\cD(F)\] from which we conclude that the sets are equal. 
\end{proof}

\begin{remark}
To highlight the analogy between the decomposition of a Boolean network and the decomposition of its dynamics, we set  
\[
\mathcal{D}(F_1)\rtimes_P \mathcal{D}(F_2) :=
\mathcal{D}(F_1)\rtimes \mathcal{D}(F_2 ^{F_1}) := \bigsqcup_{\mathcal{C}_1\in\mathcal{D}(F_1)}\mathcal{C}_1\oplus\mathcal{D}(F^{\mathcal{C}_1}_2).
\]
Theorem~\ref{thm:decomp} then states
\[\mathcal{D}(F_1\rtimes_P F_2)=\mathcal{D}(F_1)\rtimes_P \mathcal{D}(F_2),\]
which can be seen as a distributive property for dynamics of coupled networks. Combining Theorem \ref{thm:unique_decomp} and Theorem \ref{thm:decomp}, the dynamics of a Boolean network $F$, which decomposes into simple networks $F_1,\ldots,F_m$.
By a slight abuse of notation we have,
    \[\mathcal{D}(F)=\mathcal{D}(F_1)\rtimes_{P_1}\bigg(\cD(F_2)\rtimes_{P_2}\big(\cdots\rtimes_{P_{m-1}}\cD(F_m)\big)\bigg).\]
\end{remark}

\begin{example}\label{ex_conj_network}
Consider the Boolean network $F(x_1,x_2,y_1,y_2)=(x_2,x_1,x_2y_2,y_1).$ We can decompose $F$ into a network $F=F_1\rtimes F_2$ with $F_1(x_1,x_2)=(x_2,x_1)$ and $F_2(y_1,y_2)=(y_2,y_1)$. Note that $\overline{F}_2(x_1,x_2,y_1,y_2)=(x_2 y_2,y_1)$ and $\mathcal{D}(F_1)=\left\{00,11,\{01,10\}\right\}$ (note that we denote a steady state $\mathcal{C}=\{c\}$ simply by $c$). To use  Theorem~\ref{thm:decomp} we find the attractors of $F_1$ and the attractors of $F_2$ induced by each of those attractors.

\begin{itemize}
\item For $\mathcal{C}_1=\{00\},$ the corresponding network is $y(t+1)=\overline{F}_2(0,0,y(t))$.  If $y(0)=(y_1^*,y_2^*),$ then $y(1)=\overline{F}_2(0,0,y_1^*,y_2^*)=(0,y_1^*)$ and $y(2)=\overline{F}_2(0,0,0,y_1^*)=(0,0)$. Then the space of attractors for $F_2^{\mathcal{C}_1}$ is \[\mathcal{D}(F_2^{\mathcal{C}_1})=\{00\}.\]

\item For $\mathcal{C}_2=\{11\},$ the corresponding network is $y(t+1)=\overline{F}_2(1,1,y(t))$.  If $y(0)=(y_1^*,y_2^*),$ then $y(1)=\overline{F}_2(1,1,y_1^*,y_2^*)=(y_2^*,y_1^*)$  and $y(2)=\overline{F}_2(1,1,y_2^*,y_1^*)=(y_1^*,y_2^*)$. Thus, the corresponding attractor space is \[\mathcal{D}(F_2^{\mathcal{C}_2})=\{00,11, \{01,10\} \}.\] 

\item For $\mathcal{C}_3=\{01,10\},$ we define $g(t):\mathbb{N}\to X_1$ by $g(0)=\{01\},\,g(1)=\{10\},$ and $g(t+2)=g(t)$ and $F_2^{\mathcal{C}_3}$ is given by $y(t+1)=\overline{F}_2(g(t),y(t)).$ 
If $y(0)=(y_1^*,y_2^*),$ then $y(1)=\overline{F}_2(0,1,y_1^*,y_2^*)=(y_2^*,y_1^*)$,
$y(2)=\overline{F}_2(1,0,y_2^*,y_1^*)=(0,y_2^*)$,
$y(3)=\overline{F}_2(0,1,0,y_2^*)=(y_2^*,0)$,
$y(4)=\overline{F}_2(1,0,y_2^*,0)=(0,y_2^*)$.
Then, the corresponding attractor space is \[\mathcal{D}(F_2^{\mathcal{C}_3})=\{00,\{01,10\}\}.\] 
\end{itemize}

\noindent To reconstruct the space of attractors for $F$, we have
%\[\mathcal{D}(F)=
%\mathcal{D}(F_1)\times \mathcal{D}(F_2 ^{F_1})
%=(\{00\},\{11\},\{01,10\})\times \left(\mathcal{D}(F_2^{\mathcal{C}_1}),\mathcal{D}(F_2^{\mathcal{C}_2}),\mathcal{D}(F_2^{\mathcal{C}_3})\right)\]
%\[=
%\{00\}\oplus\{\{00\}\}
%+
%\{11\}\oplus\{\{00\},\{11\},\{01,10\}\}
%+
%\{01,10\}\oplus\{\{00\},\{01,10\}\]
%\[=\{\{0000\},\{1100\},\{1101,1110\},\{0100,1000\},\{0101,1010\}\}\]

\[\mathcal{D}(F)=
\mathcal{D}(F_1)\rtimes \mathcal{D}(F_2 ^{F_1})
=(00,11,\{01,10\})\rtimes \left(\mathcal{D}(F_2^{\mathcal{C}_1}),\mathcal{D}(F_2^{\mathcal{C}_2}),\mathcal{D}(F_2^{\mathcal{C}_3})\right)\]
\[=
00\oplus\{00\}
+
11\oplus\{00,11,\{01,10\}\}
+
\{01,10\}\oplus\{00,\{01,10\} \} 
\]
\[=\{0000,1100,1111,\{1101,1110\},\{0100,1000\},\{0101,1010\}\},\]

 which agrees with our calculation from Example \ref{example:small_coupled}.
\end{example}

\begin{figure}
    \centering
    \includegraphics[scale=0.4]{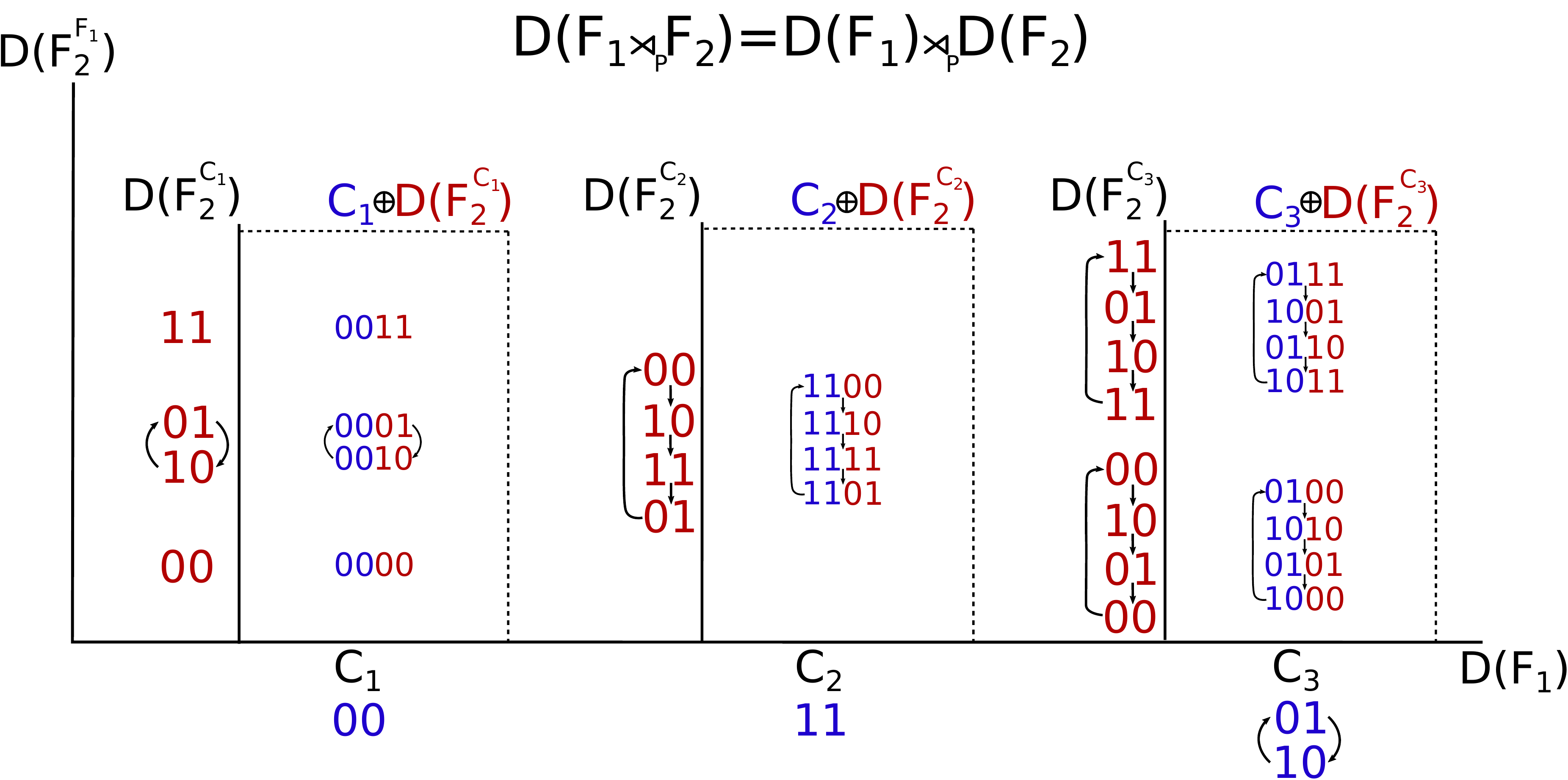}
    \caption{
    Graphical description of Theorem~\ref{thm:decomp} (applied to Example~\ref{ex_linear_network}). The dynamics of $F_1\rtimes_P F_2$ can be seen as a semi-direct product between the dynamics of $F_1$ and the dynamics of $F_2$ induced by $F_1$ via the coupling scheme $P$. 
    %The latter consists of $F_2^{\mathcal{C}_i}$ for all attractors $\mathcal{C}_i$'s of $F_1$. 
    The dynamics of $F_2$ induced by attractors of $F_1$ can vary, and Theorem~\ref{thm:decomp} shows precisely how to combine all of these attractors.
    }
    \label{fig:prod}
\end{figure}

\begin{example}\label{ex_linear_network}
Consider the linear network $F(x_1,x_2,y_1,y_2)=(x_2,x_1,x_2+y_2,y_1).$ We can decompose $F=F_1\rtimes F_2$ into simple networks $F_1(x_1,x_2)=(x_2,x_1)$ and $F_2(y_1,y_2)=(y_2,y_1)$. Note that $\overline{F}_2(x_1,x_2,y_1,y_2)=(x_2+y_2,y_1)$ and $\mathcal{D}(F_1)=\left\{00,11,\{01,10\}\right\}$.
We use Theorem~\ref{thm:decomp} as follows (see Figure~\ref{fig:prod} for a graphical description).
\\

\begin{itemize}
\item For $\mathcal{C}_1=\{00\},$ the corresponding network is $y(t+1)=\overline{F}_2(0,0,y(t))$.  If $y(0)=(y_1^*,y_2^*),$ then $y(1)=\overline{F}_2(0,0,y_1^*,y_2^*)=(y_2^*,y_1^*)$. Thus, the space of attractors for $F_2^{\mathcal{C}_1}$ is \[\mathcal{D}(F_2^{\mathcal{C}_1})=\{00,11,\{01,10\}\}.\]

\item Similarly, for $\mathcal{C}_2=\{11\},$ we find that the space of attractors for $F_2^{\mathcal{C}_2}$ is \[\mathcal{D}(F_2^{\mathcal{C}_2})=\{\{00,10,11,01\}\}\]

\item For $\mathcal{C}_3=\{01,10\},$ we define $g(t):\mathbb{N}\to X_1$ by $g(0)=\{01\},\,g(1)=\{10\},$ and $g(t+2)=g(t)$ and $F_2^{\mathcal{C}_3}$ is given by $y(t+1)=\overline{F}_2(g(t),y(t)).$ If $y(0)=(y_1^*,y_2^*)$ then $y(1)=(1+y_2^*,y_1^*),\,y(2)=(y_1^*,y_2^*+1),\,y(3)=(y_2^*,y_1^*),\,y(4)=(y_1^*,y_2^*)=y(0)$, thus in general $y(4t)=(y_1^*,y_2^*)$, $y(4t+1)=(1+y_2^*,y_1^*)$, $y(4t+2)=(y_1^*,y_2^*+1)$, and $y(4t+3)=(y_2^*,y_1^*)$ for $t>0$. It follows that there are only 2 periodic trajectories in this case: $(00,10,01,00,00,10,01,00,\ldots)$ and $(11,01,10,11,11,01,10,11,\ldots)$, which have period 4. 
The corresponding attractor space is \[\mathcal{D}(F_2^{\mathcal{C}_3})=\{\{00,10,01,00\},\{11,01,10,11\}\}.\] 
Note that the repetition of certain states is needed to obtain the correct attractors of the full network $F$.%for the period of the full attractors to be encoded. 
\end{itemize}

\noindent To reconstruct the space of attractors for $F$, we have
\[\mathcal{D}(F)=
\mathcal{D}(F_1)\rtimes \mathcal{D}(F_2 ^{F_1})
=(00,11,\{01,10\})\rtimes \left(\mathcal{D}(F_2^{\mathcal{C}_1}),\mathcal{D}(F_2^{\mathcal{C}_2}),\mathcal{D}(F_2^{\mathcal{C}_3})\right)\]
%\[
%=00\oplus\{00,11,\{01,10\}\}
%+
%11\oplus\{\{00,10,11,01\}\}
%+
%\{01,10\}\oplus\{\{00,10,01,00\},\{11,01,10,11\}\}\]
%\[=\{0000,0011,\{0001,0010\},\{1100,1110,1111,1101\},\{0100,1010,0101,1000\},\{0111,1001,0110,1011\}\}.\]

\[=\left\{\begin{matrix}
00\oplus\{00,11,\{01,10\}\}\\
11\oplus\{\{00,10,11,01\}\}\\
\{01,10\}\oplus\{\{00,10,01,00\},\{11,01,10,11\}\}
\end{matrix}\right\}=
\left\{\begin{matrix}
0000,0011,\{0001,0010\},\\
\{1100,1110,1111,1101\},\\
\{0100,1010,0101,1000\},\\
\{0111,1001,0110,1011\}
\end{matrix}\right\}.
\]

% which agrees with our calculation from Example \ref{example1}.

% From the linear systems theory we know that $\mathcal{D}(F_1)=1+C_1+C_2$ and $\mathcal{D}(F_2)=1+C_1+C_2$. However, we can also check that
% $\mathcal{D}(F)=1+C_1+C_2+3C_4$. This raises the questions if we can use Theorem~\ref{thm:decomp} to obtain $\mathcal{D}(F)$ from $\mathcal{D}(F_1)$ and $\mathcal{D}(F_2)$. This is certainly not possible without knowledge of the connection between $F_1$ and $F_2$, as exemplified by Figure~\ref{fig_linear_network} where we look at all possible non-isometric connections for this example. \claus{do we want this ending? we don't really look at this in the next section. Delete figure 6? Pose this question more openly?}
%In the following, we will investigate what we can learn about $\mathcal{D}(F)$ from $\mathcal{D}(F_1)$ and $\mathcal{D}(F_2)$ when restricting ourselves to a special, more biologically meaningful class of Boolean functions, nested canalizing functions.
\end{example}

These examples highlight how Theorem~\ref{thm:decomp} enables the computation of the dynamics of a decomposable network from the dynamics of its simple networks.

\section{Parametrizing Network Extensions}\label{sec:enumeration}
We have shown above how non-strongly connected networks can be decomposed into simple (i.e., strongly connected) networks. In this section, we consider the opposite: how and in how many ways can several (not necessarily simple) networks be combined into a larger network. For two networks, this can be thought of as determining the number of ways that the functions of one network can be extended to include inputs coming from the other network. In general, this problem can be quite difficult. However, for certain families of functions with sufficient structure we can enumerate all such extensions.%\claus{maybe rewrite slightly}  We'll look at a few cases where such an enumeration is possible.  In the first case all of our networks will consist of linear functions, and then we'll consider networks of nested canalizing functions.

\begin{definition}
Let $X = \FF_2^{n_1}, Y = \FF_2^{n_2}$. An extension of a Boolean function $f:X\to\FF_2$ by $Y$ is any Boolean function $\tilde f:X\times Y\to \FF_2$ such that $\tilde f|_X=f.$ That is, the restriction of an extension to the inputs of the original function recovers the original function.
\end{definition}

\subsection{Graphical Boolean networks}
Various families of Boolean functions feature symmetries such that a Boolean network governed by functions from such a family is uniquely determined by its wiring diagram. This includes linear functions, as well as conjunctive and disjunctive functions. We call networks governed by such a family of functions \emph{graphical}. The set of all extensions of graphical networks can be enumerated and parameterized in a straightforward manner.

% \begin{definition}
% A simple system is any Boolean network which is completely determined by its wiring diagram.
% \end{definition}

\begin{remark}\label{rem_partial_order}
The wiring diagram of a Boolean network $F$ can be represented as a square matrix with entries in $\FF_2$, the adjacency matrix. When the network is decomposable, the wiring diagram of $F$ can be represented by a lower block triangular matrix after a permutation of variables. The blocks on the diagonal represent the restrictions of $F$ to each simple network and the off-diagonal blocks represent uni-directional connections between the simple networks.  The ordering of the blocks along the diagonal corresponds to the  acyclic graph $Q$ from Definition~\ref{def:acyclic}, which defines a partial ordering of the simple networks.

Conversely, a collection of simple networks can be combined to form a larger network. The matrix representation of the wiring diagram of this larger network will consist of blocks along the diagonal representing these simple networks and blocks below the diagonal representing connections between each simple network.
\end{remark}

\begin{proposition}\label{prop:counting_graphical}
Let $X_1 = \FF_2^{n_1}, X_2 = \FF_2^{n_2}$. Given two Boolean networks $F_1: X_1\to X_1,\,F_2: X_2\to X_2$, all possible wiring diagrams of an extension of $F_2$ by $F_1$, $F = F_1 \rtimes_P F_2 : X\times Y\to X\times Y$, can be parametrized by $\FF_2^{n_2,n_1}$. There are $2^{n_1n_2}$ unique wiring diagrams for $F = F_1 \rtimes_P F_2$.
\end{proposition}
\begin{proof}
Let $W_1, W_2$ be the wiring diagrams of $F_1, F_2$. For the network $F$ to be an extension of $F_2$ by $F_1$, any $f_2\in F_2$ may depend on any $x_1\in X_1$ but the opposite is not admissible. The wiring diagram of such an extension is thus given by
    \[W=\begin{bmatrix}W_1 & 0\\ \tilde P & W_2\end{bmatrix},\]
where the submatrix $\tilde P\in \FF_2^{n_2,n_1}$  specifies the dependence of $X_2$ on $X_1$.
\end{proof}

\begin{corollary}
Let $X_1 = \FF_2^{n_1}, X_2 = \FF_2^{n_2}$. Given two graphical Boolean networks $F_1: X_1\to X_1,\,F_2: X_2\to X_2$, the number of different extensions of $F_2$ by $F_1$, $F = F_1 \rtimes_P F_2 : X\times Y\to X\times Y$, is $2^{n_1n_2}$.
\end{corollary}

\begin{example}
Recall Example~\ref{example_linear} with two graphical networks (governed by linear functions) given by $F_1(x_1,x_2)=(x_2,x_1)$ and $F_2(y_1,y_2)=(y_2,y_1)$, which were coupled together to form the larger network \[F(x_1,x_2,y_1,y_2)=(x_2,x_1,x_1+x_2+y_2,x_2+y_1).\]
In matrix form, we have
    \[F=\begin{bmatrix}F_1&0\\\tilde P&F_2\end{bmatrix},\; F_1=\begin{bmatrix}0&1\\1&0\end{bmatrix},\;F_2=\begin{bmatrix}0&1\\1&0\end{bmatrix},\;\tilde P=\begin{bmatrix}1&1\\0&1\end{bmatrix}.\]
    %\[F=\begin{bmatrix}0&1&0&0\\1&0&0&0\\1&1&0&1\\0&1&1&0\end{bmatrix}=\begin{bmatrix}F_1&0\\\tilde P&F_2\end{bmatrix}\]
The matrix $\tilde P$ is one out of a total of 16 Boolean $2\times2$ matrices, which determines the extension.
\end{example}

%We can generalize the previous results to modular networks with several modules, and start by specifying the acyclic graph that induces a partial ordering of the modules described in Remark~\ref{rem_partial_order}.

We can generalize the previous results to any network, using the acyclic graph that encodes how the simple networks are connected, described in Remark~\ref{rem_partial_order}.

% \begin{definition}
% Let $W_1,\ldots, W_m$ be the wiring diagrams of the modules of a modular Boolean network $F$. By setting $W_i \rightarrow W_j$ if there exists at least one edge from a vertex in $W_i$ to a vertex in $W_j$, we obtain a (directed) acyclic graph 
% $$Q = \{(i,j) | W_i \rightarrow W_j \},$$ which describes the connections between the modules of $F$. Note that this acyclic graph induces a partial order on the modules     (see \cite{jarrah2010dynamics}).

% %We define a partial order by setting $W_i \preccurlyeq W_j$ if there exists at least one edge from a vertex in $W_i$ to a vertex in $W_j$. The set of modules together with $\preccurlyeq$ forms a partially ordered set%, which can be expressed as a lower triangular matrix $Q$ with
% % $$q_{ij} = \begin{cases}
% %     1 & \ \text{if }W_i \preccurlyeq W_j,\\
% %     0 & \ \text{otherwise.}
% % \end{cases}$$
% \end{definition}

\begin{theorem}\label{thm:counting_graphical}
Let $X_1 = \FF_2^{n_1}, \ldots, X_m = \FF_2^{n_m}$ and let $F_1,\ldots,F_m$ with $F_i: X_i\to X_i$ be a collection of Boolean networks. All possible wiring diagrams of a network $F$ which admits a decomposition into simple networks $F_1,\ldots,F_m$ can be parametrized (i.e., are uniquely determined) by the choice of an acyclic graph $Q$ and an element in $$\prod_{\substack{(i,j)\in Q\\i\neq j}}\FF_2^{n_j,n_i}\backslash \{0\}.$$ 
\end{theorem}
\begin{proof}
Let $F$ be a network with simple networks $F_1,\ldots,F_m$ and associated acyclic graph $Q$. (Note $F$ may already be simple). Further, let $W_i$ be the wiring diagram of $F_i, i=1,\ldots, m$. For each $(i,j) \in Q$, there exists a unique non-zero matrix $\tilde P_{ij} \in \FF_2^{n_j,n_i}$, which describes the exact connections from $W_i$ to $W_j$. Note that the zero matrix would imply no connections from $W_i$ to $W_j$, contradicting $(i,j)\in Q$. 
\end{proof}

\begin{remark}
The set $\mathcal{Q}$ of all possible acyclic graphs corresponding to a Boolean network $F$ with $m\geq 1$ simple networks can be parametrized by the set of all lower triangular, Boolean matrices with diagonal entries of 1. Thus,
$|\mathcal{Q}| = 2^{m(m-1)/2}$.
\end{remark}

\begin{corollary}
The number of all possible wiring diagrams of a Boolean network with simple networks $F_1,\ldots,F_m$ in that order (i.e., with any acyclic graph $Q$ satisfying $(i,j)\not\in Q$ whenever $i>j$) is given by
\[\sum_{Q\in\mathcal{Q}}\prod_{\substack{(i,j)\in Q\\i\neq j}}(2^{n_in_j}-1).\]
Note that the product over $Q=\emptyset$ is 1. Further, note that this number can also simply be expressed as $2^M$, where 
$$M = \sum_{i=1}^m\sum_{j=i+1}^m n_in_j$$
describes the total number of possible edges between all pairs of simple networks $F_i$ and $F_j$ with $i<j$.
\end{corollary}

\begin{corollary}
Let $X_1 = \FF_2^{n_1}, \ldots, X_m = \FF_2^{n_m}$ and let $F_1,\ldots,F_m$ with $F_i: X_i\to X_i$ be a collection of graphical Boolean networks. The number of different graphical networks with simple networks $F_1,\ldots,F_m$ with acyclic graph $Q$ is given by
$$\prod_{\substack{(i,j)\in Q\\i\neq j}}(2^{n_in_j}-1),$$
and any graphical network with simple networks $F_1,\ldots,F_m$ (with any  acyclic graph  $Q$ satisfying $(i,j)\not\in Q$ whenever $i>j$) is parameterized by an element of
\[\bigsqcup_{Q\in\mathcal{Q}}\prod_{\substack{(i,j)\in Q\\i\neq j}}\FF_2^{n_j,n_i}\backslash \{0\}\cong \FF_2^M\]
where 
$$M = \sum_{i=1}^m\sum_{j=i+1}^m n_in_j.$$
\end{corollary}

% \begin{remark}
% We can generalize the previous results as follows. Let $X_1 = \FF_2^{n_1}, \ldots, X_m = \FF_2^{n_m}$. Given a collection of simple Boolean networks, $F_1: X_1\to X_1, \ldots, F_m: X_m\to X_m$, we can parameterize all extensions $F$, which have a modular decomposition into $F_1,\ldots,F_m$ by choosing an ordering $\sigma$ of the modules $F_i$, and then choosing matrices $P_{ij}$ representing the connections between $X_i$ and $X_j$ where $\sigma(i)<\sigma(j)$.
% \end{remark}
%In general, given a collection of linear systems $(A_1,X_1),\ldots, (A_m,X_m)$ we can parameterize all linear systems $A$ which have a modular decomposition into $A_1,\ldots,A_m$ by choosing an ordering $\sigma$ of the modules $A_i$, and then choosing matrices $P_{ij}$ representing the connections between $X_i$ and $X_j$ where $\sigma(i)<\sigma(j)$.

\begin{remark}
The previous results can be easily extended to networks governed by functions such that the network is uniquely determined by a signed graph, where each edge is associated with a ternary value of $-1$, $0$ (i.e., no regulation) or $1$. For example, AND-NOT networks, which are governed by functions such as $\neg x_1 \wedge  x_2$, are completely determined by a signed graph. For all such networks, the connecting matrices $\tilde P_{ij}$ in the proofs of Proposition~\ref{prop:counting_graphical} and Theorem~\ref{thm:counting_graphical} are elements of $\FF_3^{n_j,n_i}$, implying that all occurrences of $2$ ($\FF_2$) get replaced by $3$ ($\FF_3$) in the results above. %The update function at each node is determined by the product of incoming edges where a -1 corresponds to a negation of the originating node and a 1 corresponds to identity of the originating node.
%The above theory extends in a straightforward way where now the connecting matrices $\tilde P_{ij}$ in $\FF_3^{n_j,n_i}$.
\end{remark}

\begin{example}
Consider the signed graphs (in matrix form)
\[W_1=\begin{bmatrix} 0&-1\\1&1\end{bmatrix},\quad W_2=\begin{bmatrix}-1 & 1\\1 & 0\end{bmatrix}\]
with corresponding AND-NOT networks 
\[F_1=(\neg x_2,x_1\wedge x_2)\quad\text{and}\quad F_2=(\neg x_3\wedge x_4,x_3).\]
Given the connecting matrix $$\tilde P=\begin{bmatrix}0&0\\1&-1\end{bmatrix},$$ we can extend $F_2$ by $F_1$ to obtain the extended AND-NOT network%form the wiring diagram 
%\[W=\begin{bmatrix}0&-1&0&0\\1&1&0&0\\0&0&-1&1\\1&-1&1&0\end{bmatrix}\]
%describing the modular AND-NOT network 
\[F=F_1 \rtimes_P F_2 = (\neg x_2,x_1\wedge x_2,\neg x_3\wedge x_4,x_1\wedge\neg x_2\wedge x_3),\]
which naturally decomposes into $F_1$ and $F_2$. Note that $\tilde P$ is one out of a total of $3^4 = 81$ ternary $2\times 2$ matrices, each of which determines a unique extension of $F_2$ by $F_1$.
\end{example}

\subsection{Nested canalizing Boolean networks}
Thus far, we have described how to parametrize (and count) all possible extensions for networks, which are completely described by their (signed) wiring diagram. All NCFs with a single layer fall into this category of functions. For the family of all NCFs, which includes NCFs with multiple layers, this simple parametrization is no longer possible, as there are multiple possibilities how added variables can extend a general NCF. In order to understand this better, we investigate what happens when we restrict an NCF one variable at a time. To simplify notation, we let $\overbar{\{x_i\}}$ denote the set of all variables of a function except for variable $x_i$.

%The situation becomes more complicated when we consider systems of nested canalizing functions. As noted above, the reason that combinations of linear systems are so easily parameterized is due to the fact that any linear system is completely determined by its wiring diagram.  For this reason we have that say combinations of conjunctive networks are parameterized in exactly the same way as conjunctive networks are also completely determined by their wiring diagram.  Unfortunately for nested canalizing functions this is not the case.  To parameterize the extension of one nested canalizing function by another, it is not enough to specify the connections between two modules.  One needs to also specify how the NCF at each node will be extended by the newly connected nodes.  

% In order to understand how a function $f$ can be extended, we'll start by investigating what happens as we restrict an NCF one variable at a time. 
\begin{example}\label{ex_ncf_count}
Consider the following NCF with unique polynomial representation
$$f=x_1(x_2+1)\left[x_3\left[(x_4+1)x_5+1\right]+1\right]+1$$
and layer structure $(2,1,2)$. Recall that restricting an NCF to $\overbar{\{x_i\}}$ means setting $x_i$ to its non-canalizing value. We will now consider restrictions of $f$ to different $\overbar{\{x_i\}}$ to highlight how the layer structure of $f$ determines the layer structure of the restriction.
\begin{enumerate}
    \item Restricting $f$ to $\overbar{\{x_2\}}$ means setting $x_2=0$, i.e., removing the factor $(x_2+1)$ from the polynomial representation of $f$. This restriction has layer structure $(1,1,2)$ as the size of the first layer has been reduced by one.  
    \item Restricting $f$ to $\overbar{\{x_3\}}$ means setting $x_3=1$ and yields
    $$f|_{\overbar{\{x_3\}}} = x_1(x_2+1)(x_4+1)x_5+1,$$
    an NCF with a single layer. This reduction in the number of layers occurs because $x_3$ is the only variable in the second layer of $f$. Restricting to $\overbar{\{x_3\}}$ implies an annihilation of this layer and results in a fusion of the adjacent layers (layer 1 and layer 3), as all variables in these two layers share the same canalized output value. 
    \item Restricting $f$ to $\overbar{\{x_5\}}$ means setting $x_5=1$ and yields
    $$f|_{\overbar{\{x_5\}}} = x_1(x_2+1)[x_3x_4+1]+1,$$
    which is an NCF of only two layers. This occurs because the last layer of an NCF always consists of two or more variables. Removal of $x_5$ from this layer leaves only variable $x_4$, which can now be combined with the previous layer. Note that the canalizing value of $x_4$ flips in this process.
\end{enumerate}
\end{example}
\begin{remark}\label{rem:addition_rules}
Given an NCF with layer structure $(k_1,\ldots,k_r)$, removing a variable in the $i$th layer (i.e., restricting the NCF to all but this variable) results in decreasing the size of the $i$th layer by one, $k_i\mapsto k_i-1.$ The previous example highlighted the three different cases that may occur when restricting NCFs. By reversing this thought process, we can categorize the number of ways that a new variable can be added to an NCF to obtain an extension.  
\begin{enumerate}
    \item (Initial layer) A new variable can always be added as new outermost layer. That is,
        \[(k_1,\ldots,k_r)\mapsto(1,k_1,\ldots,k_r)\]
    \item (Layer addition) A new variable can always be added to any layer. That is,
        \[(k_1,\ldots,k_r)\mapsto(k_1,\ldots,k_i+1,\ldots,k_r)\]
    \item (Splitting) If $k_i\geq 2,$ then the new variable can split the $i$th layer. That is,
        \[(k_1,\ldots,k_r)\mapsto(k_1,\ldots,k_i-\ell,1,\ell,\ldots,k_r)\]
        where $1\leq\ell\leq k_i-1$, with the exception that if $i=k$ and $\ell=1$,
        \[(k_1,\ldots,k_r)\mapsto(k_1,\ldots,k_r-1,2)\]
        because the last layer of an NCF always contains at least two variables.
\end{enumerate}
 Using these rules, we can construct extensions by adding new variables one at a time.
 \end{remark}
 
\begin{lemma}\label{lem:ncf_unique_add}
Let $X=\FF_2^{n_1}$. Consider an NCF $f:X\to\FF_2$ and a set of variables $Y=\{y_1,\ldots,y_{n_2}\}$.  By fixing an ordering on $Y$, every extension $\tilde f: X\times Y\to\FF_2$ of $f$ can be realized uniquely by the sequential addition of the variables of $Y$ using the rules described in Remark~\ref{rem:addition_rules}.
\end{lemma}
\begin{proof}
The restriction of an NCF to $\overbar{\{x_i\}}$ defines a unique function, as it means setting $x_i$ to its unique non-canalizing value. Let $\tilde f: X\times Y\to \FF_2$ be an NCF extending $f$ over $Y$. By systematically restricting $\tilde f$ one variable at a time, we obtain a sequence of functions $\{f_i\}_{i=0}^{n_2}$
    \[f=f_0\to f_1\to \cdots\to f_i\to f_{i+1}\to \cdots f_{n_2}=\tilde f,\]
where each function $f_i$ is obtained by restricting $f_{i+1}$ to $\overbar{\{y_i\}}.$ As every restriction has an inverse extension, each $f_{i+1}$ can be obtained from $f_{i}$ by adding $y_{i}$ in the unique way that reverses the restriction.
\end{proof}

 \begin{figure}
    \centering
    \includegraphics[width=0.6\textwidth]{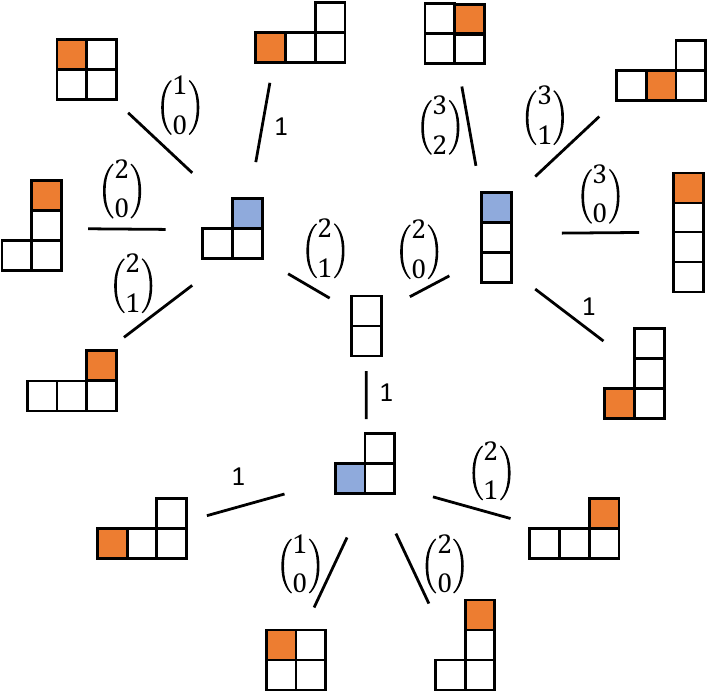}
    \caption{Graphical representation of the iterative process to find all possible nested canalizing extensions of an NCF in two variables. Each diagram represents a layer structure, each block represents a variable, each column represents a layer with the leftmost column representing the outer layer. The central diagram displays the only possible layer structure of a two-variable NCF. The colored blocks highlight the newly added variable; the first (second) added variable is shown in blue (orange) and within a layer the newly added variable is always shown on top. Edge labels show the number of extensions sharing the given layer structure. Note that for each extended layer structure, the arbitrary choice of canalizing input value for the newly added variables supplies an additional factor of 2, which the edge labels do not include.}
    \label{fig:graphical_rep}
\end{figure} 

\begin{theorem}\label{thm:ncf_enumeration}
For an NCF $f$ with layer structure $(k_1,\ldots,k_r)$, the total number of unique extensions of $f$ by one variable is given by
    \[2+2\sum_{i=1}^r(2^{k_i}-1)\]
\end{theorem}
\begin{proof}
As outlined in Remark~\ref{rem:addition_rules}, there exist three distinct ways in which a variable $y$ can be added to an NCF to obtain a new, extended NCF: 1. as a new initial (i.e., outermost) layer, 2. as an addition to an existing layer, or 3. as a splitting of an existing layer. To obtain the total number of unique extensions of $f$ by one variable, we sum up all these possibilities, separately for each layer of $f$. Note that, irrespective of the way $y$ is added, there are always two choices for its canalizing value, which is why we multiply the total number by 2 in the end.%We'll restrict our calculations to one layer at a time and we'll denote the variable to be added as $y$.  The choice of canalizing value we'll leave for the end.
% \begin{enumerate}
%     \item There exists only one way to add $y$ as new initial (i.e., outermost) layer.
%     \item For a given layer of size $k_i$, there is exactly one way to add a new variable.
%     \item 
% \end{enumerate}
%First, there is exactly one way to add $y$ as new initial (i.e., outermost) layer.

For a given layer of size $k_i$, there is exactly one way to add $y$.  If $k_i=1$, the layer cannot be split so that  there is only one way to adjust this layer. If $k_i>1,$  the layer can be split
    \[(\ldots,k_i,\ldots)\mapsto (\ldots,k_i-\ell,1,\ell,\ldots),\]
and there are $\binom{k_i}{\ell}$ ways to choose $\ell$ of the $k_i$ variables to be placed in the layer following $y$.  A proper split requires $\ell,k_i-\ell\geq 1$ so that the number of possibilities to add $y$ to a layer (by adding to the layer or splitting it) is given by
\[N_i := 1+\sum_{j=1}^{k_i-1}\binom{k_i}{j}=2^{k_i}-1.\]
Note that this formula is also valid for the special last layer of an NCF, which always contains two or more variables: Splitting the final layer into a layer of size $k_r-1$ and a new last layer of size $\ell=1$ results in this last layer to be combined with the new layer consisting only of $y$, yielding a last layer of exact size $2$, which includes $y$ (see Case 3 in Example~\ref{ex_ncf_count} and Figure~\ref{fig:graphical_rep}).

Altogether, including the sole extension which adds $y$ as initial (i.e., outermost) later and accounting for the two choices of canalizing input value of $y$, the total number of unique extensions of $f$ by one variable is thus 
\[N = 2\big(1+\sum_{i=1}^{r}N_i\big)=2+2\sum_{i=1}^{r}(2^{k_i}-1).\]
\end{proof}
By combining Lemma~\ref{lem:ncf_unique_add} and Theorem~\ref{thm:ncf_enumeration}, one can enumerate all possible extensions of an NCF by $n$ variables. Figure~\ref{fig:graphical_rep} illustrates the number of ways (up to choice of canalizing input value of the newly added variables) that two variables can be iteratively added to any NCF in two variables (which always consists of a single layer). In general, all extensions of an NCF $f$ by multiple variables can be enumerated by an iterative application of Theorem~\ref{thm:ncf_enumeration}. Starting with a two-variable NCF, there are
$$N = 8, 92, 1328, 23184, 483840, 12050112$$
such extensions by $1,\ldots,6$ variables.

%\begin{corollary}
%Let $f_1:X^m\to X^m$ and $f_2:X^n\to X^n$ be Boolean networks. Let $(x_1,\ldots,x_n)$ represent a basis of $X^n$ and let $k_\ell$ represent the number of incoming nodes for the wiring diagram of $f_2$ at $x_\ell$.  Then there are
%\[ N=\prod_{j=1}^k N_n^{k_j}\]
%networks of the form $F=f_1\rtimes f_2$.
%\end{corollary}

%%%%%%%%%%%%%%%%%%%%%%%%%%%%%%%%%%%%%%%%%%%%%
\section{Application: Control of Boolean Networks}\label{sec:control}

In this section we apply our theory of modularity to the problem of control of Boolean networks. Here, we show how the decomposition into simple networks can be used to obtain controls for each simple network, which can then be combined to obtain a control for the entire network. In this context, two types of control actions are generally considered: edge controls and node controls. For each type of control, one can consider deletions or constant expressions as defined below. The motivation for considering these control actions is that they represent the common interventions that can be implemented in practice. For instance, edge deletions can be achieved by the use of therapeutic drugs that target specific gene interactions while node deletions represent the blocking of effects of products of genes associated to these nodes; see~\cite{choi2012attractor,wooten2021mathematical}.

Once the simple networks have been identified, different methods for phenotype control (that is, control of the attractor space) can be used to identify controls in these networks. There include methods using Stable Motifs~\cite{zanudo2015cell}, Feedback Vertex Sets~\cite{zanudo2017structure}, as well as algebraic methods~\cite{Murrugarra:2016ul,Murrugarra:2015ud,Sordo-Vieira:2019uf}. For our examples below, we will use ~\cite{zanudo2015cell}, ~\cite{Murrugarra:2016ul}, and ~\cite{zanudo2017structure} to find controls in the simple networks. 

A Boolean network $F=(f_1,\ldots,f_n):\FF_2^n\rightarrow \FF_2^n$ with \emph{control} is %given by a function
a Boolean network $\mathcal{ F}:\mathbb{F}_2^n\times U\rightarrow \mathbb{F}_2^n$,
where $U$ is a set that denotes all possible controls, defined below. The case of no control coincides with the original Boolean network, that is, $\mathcal{F}(x,0)=F(x)$. Given a control $u\in U$, the dynamics are given by $x(t+1)=\mathcal{F}(x(t),u)$. See~\cite{Murrugarra:2016ul} for additional details and examples of how to encode edges and nodes in a Boolean network.

% %A Boolean network with \emph{control} is given by a function $\mathcal{ F}:\mathbb{F}_2^n\times U\rightarrow \mathbb{F}_2^n$,
% where $U$ is a set that denotes all possible control inputs. Given a control $u\in U$, the dynamics are given by $x(t+1)=\mathcal{F}(x(t),u)$.
% %%%%%%%%%%%%%%%%%%%%%%%%%%%%%%%%%%%%%%%%%%%%%
% %%%%%%%%%%%%%%%%%%%%%%%%%%%%%%%%%%%%%%%%%%%%%

% We consider a Boolean network $F=(f_1,\ldots,f_n):F_2^n\rightarrow \mathbb{F}_2^n$ and show how to encode edge and node controls by $\mathcal{F}:\mathbb{F}_2^n\times U \rightarrow \mathbb{F}_2^n$, such that $\mathcal{F}(x,0)=F(x)$. That is, the case of no control coincides with the original BN. We remark that these control types can be combined, but for clarity we present them separately. See~\cite{Murrugarra:2016ul} for additional details and examples of how to encode edges and nodes in a Boolean network.

\begin{definition}[Edge Control, \cite{Murrugarra:2016ul}]\label{def:edge_del}
Consider the edge $x_i\rightarrow x_j$ in the wiring diagram ${W}$. 
The function
\begin{equation}
\label{edge_del_def}
\mathcal{F}_j(x,u_{i,j}) := f_j(x_1,\dots,(u_{i,j}+1)x_i+u_{ij}a_i,\dots,x_n),
\end{equation}
where $a_i$ is a constant in $\FF_2$, encodes the control of the edge $x_i\rightarrow x_j$, since for each possible value of $u_{i,j}\in \mathbb{F}_2$ we have the following control settings: 
\begin{itemize}
  \item If $u_{i,j}=0$, $\mathcal{F}_j(x,0) = f_j(x_1,\dots,x_i,\dots,x_n)$. That is, the control is not active.
  \item If $u_{i,j}=1$, $\mathcal{F}_j(x,1) = f_j(x_1,\dots,x_i=a_i,\dots,x_n)$. In this case, the control is active, and the action represents the removal of the edge $x_i\rightarrow x_j$ when $a_i=0$, and the constant expression of the edge if $a_i=1$. We use  $x_i\xrightarrow[]{a_i} x_j$ to denote that the control is active.
\end{itemize}

\end{definition}
This definition can easily be extended for the control of many edges, so that we obtain $\mathcal{F}:\mathbb{F}_2^n\times \mathbb{F}_2^e \rightarrow \mathbb{F}_2^n$, where $e$ is the number of edges in the wiring diagram. Each coordinate, $u_{i,j}$, of $u$ in $\mathcal{F}(x,u)$ encodes the control of an edge $x_i\rightarrow x_j$.
\begin{definition}[Node Control, \cite{Murrugarra:2016ul}]\label{def:node_del}
Consider the node $x_i$ in the wiring diagram ${W}$. The function
\begin{equation}
\label{node_del_def}
\mathcal{F}_j(x,u^{-}_i,u^{+}_i) := (u^{-}_i+u^{+}_i+1)f_j(x) + u^{+}_i
\end{equation}
encodes the control (knock-out or constant expression) of the node $x_i$, since for each possible value of $(u^{-}_i,u^{+}_i)\in \mathbb{F}_2^2$ we have the following control settings:  
\begin{itemize}
  \item For $u^{-}_i=0, u^{+}_i=0$, $\mathcal{F}_j(x,0,0) = f_j(x)$. That is, the control is not active.
  \item For $u^{-}_i=1, u^{+}_i=0$, $\mathcal{F}_j(x,1,0)  = 0$. This action represents the knock out of the node $x_j$.
  \item For $u^{-}_i=0, u^{+}_i=1$, $\mathcal{F}_j(x,0,1)  = 1$. This action represents the constant expression of the node $x_j$.
  \item For $u^{-}_i=1, u^{+}_i=1$, $\mathcal{F}_j(x,1,1) =  f_j(x_{t_1},\dots,x_{t_m})+1$. This action changes the Boolean function to its negative value.
\end{itemize}
\end{definition}

From Theorem~\ref{thm:decomp} the attractors of $F$ are of the form $\mathcal{C} = \mathcal{C}_1\oplus\mathcal{C}_2$ where $\mathcal{C}_1\in\mathcal{D}(F_1)$ and $\mathcal{C}_2\in\mathcal{D}(F^{\mathcal{C}_1}_2)$. Thus, we have the following theorem.

\begin{definition}
We say that a control $\mu$ stabilizes a network $F$ at an attractor $C$ when the resulting network after applying $\mu$ into $F$ (denoted here as $F^\mu$) has $C$ as its only attractor.
\end{definition}

\begin{theorem}\label{thm:control0}
Given a decomposable network $F = F_1\rtimes_P F_2$. Let $\mathcal{C} = \mathcal{C}_1\oplus\mathcal{C}_2$ be an attractor of $F$, where $\mathcal{C}_1\in\mathcal{D}(F_1)$ and $\mathcal{C}_2\in\mathcal{D}(F^{\mathcal{C}_1}_2)$ and either $\mathcal{C}_1$ or $\mathcal{C}_2$ is a steady state. If $\mu_1$ is a control that stabilizes $F_1$ in $\mathcal{C}_1$ and $\mu_2$ is a control that stabilizes $F^{\mathcal{C}_1}_2$ in $\mathcal{C}_2$, then $\mu=(\mu_1,\mu_2)$ is a control that stabilizes $F$ in $\mathcal{C}$.
\end{theorem}
\begin{proof}
Let $F^{\mu_1}_1$ be the resulting network after applying the control $\mu_1$. Thus, the dynamics of $F^{\mu_1}_1$ is $\mathcal{C}_1$, that is $\mathcal{D}(F^{\mu_1}_1) =\mathcal{C}_1$. Similarly, the %induced
dynamics of $F^{ \mathcal{C}_1, \mu_2}_2$ %by $F^{\mu_1}_1$ 
is $\mathcal{C}_2$. That is, $\mathcal{D}(F^{\mathcal{C}_1, \mu_2 }_2) = \mathcal{C}_2$. Then,
\[
F^\mu = (F_1\rtimes_P F_2)^\mu = F_1^\mu\rtimes_P F_2^\mu = F_1^{\mu_1}\rtimes_P F_2^{\mu_2}
\]
Thus,
\[
\mathcal{D}(F^\mu) = \mathcal{D}(F_1^{\mu_1}\rtimes_P F_2^{\mu_2}) = \bigsqcup_{\mathcal{C}'\in\mathcal{D}(F_1^{\mu_1})}\mathcal{C}'\oplus\mathcal{D}(F^{\mathcal{C}',\mu_2}_2)
= \mathcal{C}_1\oplus\mathcal{D}(F^{\mathcal{C}_1,\mu_2}_2) = \mathcal{C}_1 \oplus \mathcal{C}_2.
\]
It follows that there is only one attractor of $F^\mu$ and that attractor is $\mathcal{C}_1\oplus\mathcal{C}_2$.
Thus, $F$ is stabilized by $\mu=(\mu_1, \mu_2)$ and we have $\mathcal{D}(F^{\mu}) = \mathcal{C}$.
\end{proof}

\begin{remark}\label{rem:control0}
In the proof of Theorem~\ref{thm:control0} we used the fact the product of a steady state and a cycle (or vice versa) will result in only one attractor for the composed network. The former is not always true in general because multiplying two attractors (of length greater than 1) might result in several attractors for the composed network due to the attractors starting at different states. 
Nevertheless, Theorem~\ref{thm:control0} is still valid if we use the following definition of stabilization for non-autonomous networks which will guarantee that $\mathcal{C}_1$ and $\mathcal{C}_2$ can be combined in a unique way resulting in a unique attractor of the whole network. %don't erase this sentence, just comment if needed
\end{remark}
\begin{definition} %don't erase, just comment if needed
Consider a non-autonomous network given by $y(t+1)=\overline{F}_2(g(t),y(t),u)$, where $g(t)$ is a trajectory representation of an attractor $\mathcal{C}_1$ of an upstream network. We say that a control $\mu_2$ stabilizes this network, $F_2^{\mathcal{C}_1}$,  at an attractor $\mathcal{C}_2$ when the resulting network after applying $\mu_2$ (denoted here as $F_2^{\mathcal{C}_1,\mu_2}$) has $\mathcal{C}_2$ as its unique attractor. 
For non-autonomous networks the definition of unique attractor requires that $\left(g(t),y(t)\right)_{t=0}^\infty$ has a unique periodic trajectory up to shifting of $t$ (which is automatically satisfied if $\mathcal{C}_1$ or $\mathcal{C}_2$ is a steady state).
\end{definition}

\begin{example}\label{ex:control0}
Consider the network $F(x_1,x_2,x_3,x_4)=(x_2,x_1,x_2x_4,x_3)$ which can be decomposed into $F=F_1\rtimes F_2$, with $F_1(x_1,x_2)=(x_2,x_1)$ and $F_2(x_3,x_4)=(x_4,x_3)$ and 
suppose that we want to stabilize $F$ in 1111, which is an attractor (but the only one) of $F$. Note that $\overline{F}_2(x_1,x_2,x_3,x_4)=(x_2x_4,x_3)$ and $\mathcal{D}(F_1)=\left\{00,11,\{01,10\}\right\}$. Let $\mathcal{C}_1=\{11\}\in\mathcal{D}(F_1)$.
\begin{itemize}
\item The edge control $\mu_1:x_1\xrightarrow[]{1} x_2$ (that is, the control that constantly expresses the edge from $x_1$ to $x_2$) stabilizes $F_1$ at
$\mathcal{C}_1=\{11\}$. Note that the space of attractors for $F_2^{\mathcal{C}_1}$ is $\mathcal{D}(F_2^{\mathcal{C}_1})=\{00,11,\{01,10\}\}$.

\item The edge control $\mu_2:x_4\xrightarrow[]{1} x_3$ (that is, the control that constantly expresses the edge from $x_4$ to $x_3$) stabilizes $F_2^{\mathcal{C}_1}$ at $\mathcal{C}_2=\{11\}\in\mathcal{D}(F_2^{\mathcal{C}_1})$.

\item Now, the control $\mu=(\mu_1,\mu_2) = (x_1\xrightarrow[]{1} x_2, x_4\xrightarrow[]{1} x_3)$ stabilizes $F$ at $\mathcal{C} = \mathcal{C}_1\oplus\mathcal{C}_2 = \{1111\}$.
\end{itemize}
\end{example}

We note that Theorem~\ref{thm:control0} provides a way to identify controls using the modular structure of networks to stabilize the system into existing attractors. More generally, we have the following theorem.

\begin{theorem}\label{thm:control1}
Given a decomposable network $F = F_1\rtimes_P F_2$. If $\mu_1$ is a control that stabilizes $F_1$ in $\mathcal{C}_1$ (whether $\mathcal{C}_1$ is an existing attractor or a new one) and $\mu_2$ is a control that stabilizes $F^{\mathcal{C}_1}_2$ in $\mathcal{C}_2$ (whether $\mathcal{C}_2$ is an existing attractor or a new one), then $\mu = (\mu_1, \mu_2)$ is a control that stabilizes $F$ 
in $\mathcal{C} = \mathcal{C}_1\oplus\mathcal{C}_2$ provided that either $\mathcal{C}_1$ or $\mathcal{C}_2$ is a steady state.
\end{theorem}
\begin{proof}
The proof is identical to the proof of Theorem~\ref{thm:control0}.
\end{proof}

\begin{example}
Consider again the Boolean network from Example~\ref{ex:control0}, $F(x_1,x_2,x_3,x_4)=(x_2,x_1,x_2x_4,x_3)=F_1\rtimes F_2$ with $F_1(x_1,x_2)=(x_2,x_1)$ and $F_2(x_3,x_4)=(x_4,x_3)$ and 
suppose that we want to stabilize $F$ in 0110, which is not an attractor.
Note that $\overline{F}_2(x_1,x_2,x_3,x_4)=(x_2x_4,x_3)$ and $\mathcal{D}(F_1)=\left\{00,11,\{01,10\}\right\}$.
\begin{itemize}
\item Consider the control $\mu_1:(x_1\xrightarrow[]{1} x_2,x_2\xrightarrow[]{0} x_1)$. That is, the control is the combined action of setting the input from $x_1$ to $x_2$ to 1 and the input from $x_2$ to $x_1$ to 0. The control $\mu_1$ stabilizes $F_1$ at
$\{01\}$, which is not an original attractor of $F_1$. Let $\mathcal{C}_1=\{01\}\in\mathcal{D}(F_1^{\mu_1})$.
Note that the space of attractors for $F_2^{\mathcal{C}_1}$ is $\mathcal{D}(F_2^{\mathcal{C}_1})=\{00,11,\{01,10\}\}$.

\item Now consider the control $\mu_2:(x_4\xrightarrow[]{1} x_3,x_3\xrightarrow[]{0} x_4)$. That is, the control is the combined action of setting the input from $x_4$ to $x_3$ to 1 and the input from $x_3$ to $x_4$ to 0). This control stabilizes $F_2^{\mathcal{C}_1}$ at $\mathcal{C}_2=\{10\}\in\mathcal{D}(F_2^{\mathcal{C}_1})$, which is not an original attractor of $F_2^{\mathcal{C}_1}$.

\item Finally, the control $\mu=(\mu_1, \mu_2)$ stabilizes $F$ at $\mathcal{C} = \mathcal{C}_1\oplus\mathcal{C}_2 = \{0110\}$. Note that $\mathcal{C}$ is a new attractor for $F$.
\end{itemize}
\end{example}
Next we will decompose a published Boolean network model into its simple networks (corresponding to the biological modules), identify controls for each simple network and show how this yields a set of controls that directs the total system into one of its attractors.
\begin{example}
%Consider the multicellular Boolean network model (Figure~\ref{fig:pcc_model}), published in~\cite{plaugher2021modeling}, which describes the interactions of pancreatic cancer cells (purple nodes), pancreatic stellate cells (blue nodes), and their connecting cytokines (yellow nodes). This network has 69 nodes and 114 edges. We identified the nontrivial modules of this network, see Figure~\ref{fig:pcc_model}. The most upstream module is the box in amber. The green module is downstream the amber module, and the grey module is downstream the green module. The upstream module (in amber) has four attractors: two steady states and two 3-cycles. Since the upstream module consists of two triangles joined by the node TGFb1, it is enough to control this node to stabilize this module into any of its attractors. Next, we use the middle module to identify the next control targets. We identified two control nodes in the green module which are $RAS$ in the pancreatic cell (in purple) and also $RAS$ in the stellate cell (in blue). Finally, after applying the controls from the upstream modules, the last module becomes trivial and does not require additional controls. Thus, we identified the three nodes within the modules that can be used to control the entire system. 
Consider the multicellular Boolean network model, published in~\cite{plaugher2021modeling}, which studies the microenvironment of pancreatic cancer cells by modeling the interactions of pancreatic cancer cells (PCCs), pancreatic stellate cells (PSCs), and their connecting cytokines. This network has 69 nodes and 114 edges, and possesses three non-trivial simple networks, corresponding to biological modules (Figure~\ref{fig:pcc_model}a). The directed acyclic graph, which describes the connections between the modules, is shown in Figure~\ref{fig:pcc_model}b. 

% We note that the first two upstream modules contain autocrine loops involving the signaling of cytokines $TGFb1$, $EGF$, $bFGF$, and $VEGF$ which are known to be key modulators for the interaction between PCCs and PSCs~\cite{aguilar2020generalizable}. 
Using these simple networks, we can identify controls that will stabilize the network in a desired attractor. In this example, we use the phenotype apoptosis as our desired attractor. The upstream module (highlighted in amber) has four attractors: two steady states and two 3-cycles. Since this module consists of two feedback loops of length 3 joined by the node $TGFb1$, it is enough to control this node (by \cite{zanudo2017structure}) to stabilize this module into any of its attractors. Next, we move to the middle module (highlighted in green) to identify the next control targets. Using the methods from \cite{zanudo2015cell} or \cite{Murrugarra:2016ul}, we can identify a minimal set of two nodes, which need to be controlled to stabilize this module: $RAS$ in the pancreatic cell (in purple) and also $RAS$ in the stellate cell (in blue). Finally, after applying the controls from the upstream modules, the nodes in the last module (highlighted in gray) all become constant
and therefore they do not require additional controls. Thus, we identified, using the modular structure of the network, three nodes, which suffice to control the entire network. 
\begin{figure}
    \centering
    \includegraphics[width=\textwidth]{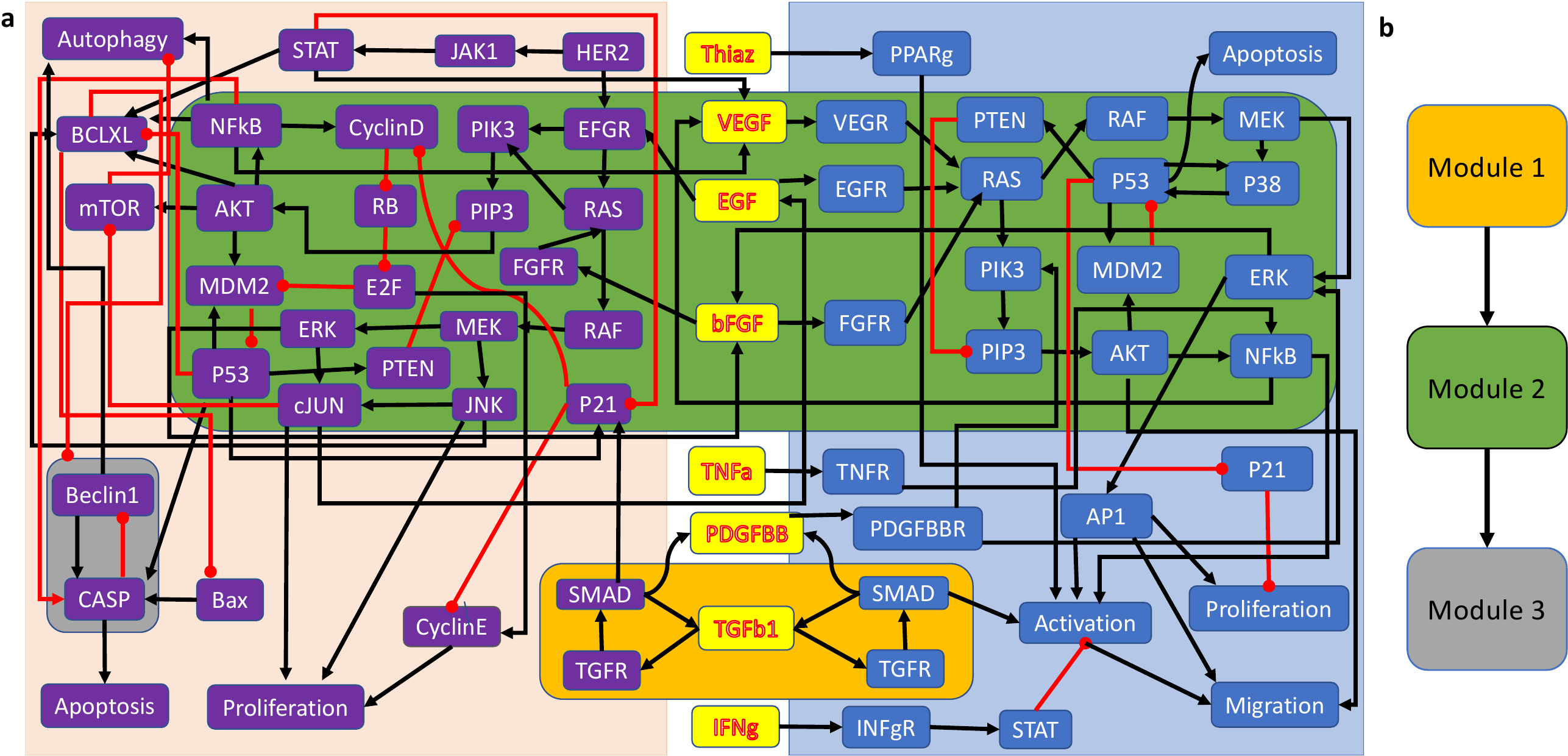}
    \caption{(a) Wiring diagram of a multicellular Boolean network model, published in \cite{plaugher2021modeling}, which  describes the interactions of pancreatic cancer cells (purple nodes), pancreatic stellate cells (blue nodes), and their connecting cytokines (yellow nodes). The non-trivial simple networks, corresponding to biological modules, are
    highlighted by amber, green, and gray boxes. (b) Directed acyclic graph describing the connections between the non-trivial simple networks.}
    \label{fig:pcc_model}
\end{figure}

\end{example}

\begin{remark}
This example highlights how one can use the decomposition into simple networks to efficiently identify a set of three controls that suffices to stabilize the whole network in a desired attractor. When disregarding the decomposition and identifying controls for the whole network instead, we also identify the same minimal set of controls. Note that this may not always be the case. In other words, the algorithm for phenotype control described here is not guaranteed to find a control of minimal size. 
\end{remark}

\section{Discussion}

Decomposition theories are a common theme across mathematics, such as group theory, low-dimensional topology, differential geometry, and other fields. The decomposition results in this paper follow this theme, applied to certain types of dynamical systems. The results are motivated by strong biological evidence that such decompositions should exist in biological systems that can be modeled by these dynamical systems. Thus, this work can be viewed as an example of biology inspiring new mathematics. The feedback to biology of our work is to provide a rigorously specified framework that can be used to make more specific the concept of a module and the modular structure of biological systems, something that is currently lacking. Furthermore, it focuses on their dynamic nature, also a feature that is currently largely absent. We have chosen the modeling framework of Boolean networks, commonly used to model gene regulatory networks. We further restrict ourselves to Boolean networks that are composed of fanout-free or linear Boolean functions, classes of functions that are restrictive enough to have interesting mathematical properties, but general enough to capture almost all biological applications of Boolean networks to gene regulatory networks. 

The core of our contributions is the definition of a module and a composition of modules into larger networks. The key result is a theorem that asserts that for a general fanout-free Boolean network, there exists a decomposition into modules, and the collection of modules is unique. Conversely, we show how one can build general networks out of modules and count the number of ways in which this can be done. A module is defined by the structural requirement that its wiring diagram is strongly connected. This condition has also been studied in the context of biochemical reaction network models consisting of systems of ordinary differential equations. There, networks with this condition are called ``weakly reversible" \cite{anderson2010weaklyreversible}, and are shown to have special properties. This suggests that it might be of interest to carry out a decomposition program similar to ours for biochemical reaction networks.

As mentioned in the introduction, there is a substantial literature on why and how modularity of biological systems could have evolved, and how it could be advantageous to an organism~\cite{borriello2021basis,manicka2021biological}. We provide one possible reason by looking at the concept of phenotype control. Here, a phenotype of a gene regulatory network is defined to be an attractor of the Boolean network model. Control inputs consist of certain modifications of the nodes and edges of the network or, alternatively, the Boolean functions on the network nodes. The aim is to find control inputs that change the attractor that a given collection of model states reaches. For instance, it is desirable to use a drug to change the cancer phenotype of a cell to a noncancer phenotype under ``normal" conditions. We show in this paper how the modular structure of a network can be used to assemble a control strategy one module at a time and to capitalize on the effect of controls in one module on downstream modules.

We note that the concept of simple networks presented in Section~\ref{subsec:simple_networks} is related to other concepts that have been used for the analysis and control of Boolean networks. Mainly, stable motifs~\cite{zanudo2015cell} are strongly connected components in an expanded graph of a Boolean network that is used for the identification of attractors and controls. For special classes of networks, the concept of simple network coincides with the one of stable motifs. For instance, when all the regulatory functions are linear functions, then the stable motifs are the same as the simple networks. We also note that the concept of simple networks presented here is related to the concept of basic blocks described in~\cite{paul2018decomposition} that has been used for the identification of attractors and controls. The difference between a simple network and a basic block is that the latter is a strongly connected component together with all upstream variables, while a simple network is just the strongly connected component.
% ``A Decomposition-based Approach towards the Control of Boolean Networks"

Several open problems remain to be tackled going forward. Most importantly, our decomposition theory needs to be related to actual regulatory networks and their structure. An important mathematical problem is to understand the dynamics of simple fanout-free or nested canalizing Boolean networks. This has been done in \cite{jarrah2010dynamics} for the special case of conjunctive networks, those constructed using only the Boolean AND function.

\section*{Acknowledgements}
The authors thank the Banff International Research Station for support through its Focused Research Group program during the week of May 29, 2022 (22frg001), which was of great help in completing this paper. M.W. thanks the University of Florida College of Medicine for travel support. The authors thank Elena Dimitrova for many fruitful initial discussions.

\bibliographystyle{siamplain}
\bibliography{main}

\end{document}

%% file: main_shared.tex
\usepackage{lipsum}
\usepackage{amsfonts}
\usepackage{graphicx}
\usepackage{epstopdf}
\usepackage{algorithmic}
\ifpdf
  \DeclareGraphicsExtensions{.eps,.pdf,.png,.jpg}
\else
  \DeclareGraphicsExtensions{.eps}
\fi

\usepackage{amssymb}
\usepackage{comment}
\usepackage{enumerate}
\usepackage{pgfplots}
\pgfplotsset{compat=1.12}
\usetikzlibrary{fit,shapes.misc,arrows,positioning,calc,decorations,decorations.pathreplacing}
\usetikzlibrary{arrows.meta}
\usetikzlibrary{shapes.geometric}

\newcounter{comments}
\definecolor{Red}{rgb}{0.8,0,0}
\definecolor{Red}{rgb}{0.8,0,0}
\definecolor{Green}{rgb}{0.2,0.6,0.2}
\definecolor{Blue}{rgb}{0,0,0.8}

\newcommand{\cC}{\mathcal{C}}
\newcommand{\cD}{\mathcal{D}}
\newcommand{\F}{\mathbb F_2}

\newcommand{\pa}{\text{pr}_1}
\newcommand{\pb}{\text{pr}_2}

\DeclareMathOperator{\vecx}{\mathbf{x}}
\DeclareMathOperator{\FF}{\mathbb{F}}

\def\longto{\longrightarrow}

\newcommand{\overbar}[1]{\mkern 1.5mu\overline{\mkern-1.5mu#1\mkern-1.5mu}\mkern 1.5mu}

\usepackage{enumitem}
\setlist[enumerate]{leftmargin=.5in}
\setlist[itemize]{leftmargin=.5in}

\newsiamremark{remark}{Remark}
\newsiamremark{hypothesis}{Hypothesis}
\newsiamremark{example}{Example}
\crefname{hypothesis}{Hypothesis}{Hypotheses}
\newsiamthm{claim}{Claim}

\headers{Decomposition of Boolean networks}{C. Kadelka, R. Laubenbacher, D. Murrugarra, A. Veliz-Cuba, and M. Wheeler}

\title{Decomposition of Boolean networks: An approach to modularity of biological systems\thanks{Currently under peer review.\funding{C.K. was partially supported by the Simons Foundation grant 712537. R.L. was partially supported by the following grants: NIH grants R011AI135128-01, R01GM127909-01, and U01EB024501-01; and NSF grant CBET-1750183.
D.M. was partially supported by a Collaboration grant (850896) from the Simons Foundation.
A.VC. was partially supported by the Simons Foundation grant 516088. M.W. was partially supported under the Keck Foundation Medical Research Grant 994413.
}}}

\author{Claus Kadelka\thanks{Department of Mathematics, Iowa State University, Ames, IA 
  (\email{ckadelka@iastate.edu}, \url{https://faculty.sites.iastate.edu/ckadelka/}).}
\and Reinhard Laubenbacher\thanks{Department of Medicine, University of Florida, Gainesville, FL 
  (\email{reinhard.laubenbacher@medicine.ufl.edu}, \url{https://directory.ufhealth.org/laubenbacher-reinhard}).}
\and David Murrugarra\thanks{Department of Applied Mathematics, University of Kentucky, Lexington, KY 
  (\email{murrugarra@uky.edu}, \url{https://sites.google.com/view/murrugarra/home}).}
\and Alan Veliz-Cuba\thanks{Department of Mathematics, University of Dayton, Dayton, OH 
  (\email{avelizcuba1@udayton.edu}).}
  \and Matthew Wheeler\thanks{Department of Medicine, University of Florida, Gainesville, FL (\email{mwheeler1@ufl.edu}).}}

\usepackage{amsopn}